\documentclass[11pt,amsfonts]{amsart}

\usepackage{fullpage}

\newtheorem{theorem}{Theorem}[section]

\newtheorem{corollary}[theorem]{Corollary}

\newtheorem{lemma}[theorem]{Lemma}

\newtheorem{proposition}[theorem]{Proposition}

\newtheorem{remark}[theorem]{Remark}

\newtheorem{example}[theorem]{Example}

\def\endproof{\qed \medskip}
\def\blacksquare{\hbox to .60em{\vrule width .60em height .60em}}

\begin{document}

\title[]{Unique continuation results for Ricci curvature and 
Applications}

\author[]{Michael T. Anderson}

\thanks{Partially supported by NSF Grant DMS 0305865 and 0604735 \\
MSC Classification: 58J32, 58J60, 53C21. Keywords: Einstein metrics, 
unique continuation }

\abstract{Unique continuation results are proved for metrics with prescribed 
Ricci curvature in the setting of bounded metrics on compact manifolds with 
boundary, and in the setting of complete conformally compact metrics on such 
manifolds. Related to this issue, an isometry extension property is proved: 
continuous groups of isometries at conformal infinity extend into the bulk of 
any complete conformally compact Einstein metric. Relations of this property 
with the invariance of the Gauss-Codazzi constraint equations under deformations 
are also discussed.}
\endabstract

\maketitle

\setcounter{section}{0}

\section{Introduction.}
\setcounter{equation}{0}

 In this paper, we study certain issues related to the boundary 
behavior of metrics with prescribed Ricci curvature. Let $M$ be a 
compact $(n+1)$-dimensional manifold with compact non-empty boundary 
$\partial M$. We consider two possible classes of Riemannian metrics 
$g$ on $M$. First, $g$ may extend smoothly to a Riemannian metric on 
the closure $\bar M = M \cup  \partial M$, thus inducing a Riemannian 
metric $\gamma  = g|_{\partial M}$ on $\partial M$. Second, $g$ may be 
a complete metric on $M$, so that $\partial M$ is ``at infinity''. In 
this case, we assume that $g$ is conformally compact, i.e.~there exists 
a defining function $\rho$ for $\partial M$ in $M$ such that the 
conformally equivalent metric
\begin{equation} \label{e1.1}
\widetilde g = \rho^{2}g 
\end{equation}
extends at least $C^{2}$ to $\partial M$. The defining function $\rho$ is 
unique only up to multiplication by positive functions; hence only the 
conformal class $[\gamma]$ of the associated boundary metric $\gamma  = 
\bar g|_{\partial M}$ is determined by $(M, g)$. 

 The issue of boundary regularity of Riemannian metrics $g$ with 
controlled Ricci curvature has been addressed recently in several 
papers. Thus, [4] proves boundary regularity for bounded metrics $g$ 
on $M$ with controlled Ricci curvature, assuming control on the 
boundary metric $\gamma$ and the mean curvature of $\partial M$ in $M$. 
In [16], boundary regularity is proved for conformally compact Einstein 
metrics with smooth conformal infinity; this was previously proved by 
different methods in dimension 4 in [3], cf.~also [5].

\medskip

 One purpose of this paper is to prove a unique 
continuation property at the boundary $\partial M$ for bounded metrics or for 
conformally compact metrics. We first state a version of the result for Einstein 
metrics on bounded domains.
\begin{theorem} \label{t 1.1.}
  Let $(M, g)$ be a $C^{2,\alpha}$ metric on a compact manifold with boundary 
$M$, with induced metric $\gamma  = g|_{\partial M}$, and let $A$ be the $
2^{\rm nd}$ fundamental form of $\partial M$ in $M$. Suppose the Ricci 
curvature $Ric_{g}$ satisfies
\begin{equation} \label{e1.2}
Ric_{g} = \lambda g, 
\end{equation}
where $\lambda$ is a fixed constant.

 Then $(M, g)$ is uniquely determined up to local isometry by the Cauchy data 
$(\gamma, A)$ on an arbitrary open set $U$ of $\partial M$. In particular, the 
topology of $M$ and $\partial M$ are uniquely determined up to covering spaces, 
and the global Cauchy data $(\gamma, A)$ on $\partial M$ are determined by their 
values on $U$.
\end{theorem}

  Similar results hold for metrics which satisfy other covariant equations involving 
the metric to $2^{\rm nd}$ order, for example the Einstein equations coupled to 
other fields; see Proposition 3.6.

\medskip

  For conformally compact metrics, the $2^{\rm nd}$ fundamental form $A$ of the 
compactified metric $\bar g$ in (1.1) is umbilic, and completely determined by 
the defining function $\rho$. In fact, for conformally compact Einstein metrics, 
the higher order Lie derivatives ${\mathcal L}_{N}^{(k)}\bar g$ at $\partial M$, where 
$N$ is the unit vector in the direction $\bar \nabla \rho$, are determined by the 
conformal infinity $[\gamma]$ and $\rho$ up to order $k < n$. Supposing $\rho$ is a 
geodesic defining function, so that $||\bar \nabla \rho|| = 1$, let 
\begin{equation}\label{e1.3}
g_{(n)} = \tfrac{1}{n!}{\mathcal L}_{N}^{(n)}\bar g.
\end{equation}
More precisely, $g_{(n)}$ is the $n^{\rm th}$ term in the Fefferman-Graham expansion 
of the metric $g$; this is given by (1.3) when $n$ is odd, and in a similar way when $n$ 
is even, cf. [18] and \S 4 below. The term $g_{(n)}$ is the natural analogue of $A$ for 
conformally compact Einstein metrics.

\begin{theorem} \label{t 1.2.}
  Let $g$ be a $C^{2}$ conformally compact Einstein metric on a compact manifold $M$ with 
$C^{\infty}$ smooth conformal infinity $[\gamma]$, normalized so that
\begin{equation} \label{e1.4}
Ric_{g} = -ng, 
\end{equation}
 Then the Cauchy data $(\gamma, g_{(n)})$ restricted to any open set $U$ of 
$\partial M$ uniquely determine $(M, g)$ up to local isometry and determine 
$(\gamma, g_{(n)})$ globally on $\partial M$.
\end{theorem}

  The recent boundary regularity result of Chru\'sciel et al.,~[16], implies that $(M, g)$ 
is $C^{\infty}$ polyhomogeneous conformally compact, so that the hypotheses of Theorem 1.2 
imply the term $g_{(n)}$ is well-defined on $\partial M$. A more general version of 
Theorem 1.2, (without the smoothness assumption on $[\gamma]$), is proved in \S 4, 
cf. Theorem 4.1. For conformally compact metrics coupled to other fields, 
see Remark 4.5. 

\medskip

  Of course neither Theorem 1.1 or 1.2 hold when just the boundary metric 
$\gamma$ on $U \subset \partial M$ is fixed. For example, in the context of 
Theorem 1.2, by [20] and [16], given any $C^{\infty}$ smooth boundary metric $\gamma$ 
sufficiently close to the round metric on $S^{n}$, there is a smooth (in the polyhomogeneous 
sense) conformally compact Einstein metric on the $(n+1)$-ball $B^{n+1}$, close to the 
Poincar\'e metric. Hence, the behavior of $\gamma$ in $U$ is independent of its 
behavior on the complement of $U$ in $\partial M$. 

  Theorems 1.1 and 1.2 have been phrased in the context of ``global'' Einstein 
metrics, defined on compact manifolds with compact boundary. However, the proofs 
are local, and these results hold for metrics defined on an open manifold with 
boundary. From this perspective, the data $(\gamma, A)$ or $(\gamma, g_{(n)})$ 
on $U$ determine whether Einstein metric $g$ has a global extension to an 
Einstein metric on a compact manifold with boundary, (or conformally compact Einstein 
metric), and how smooth that extension is at the global boundary. 

\medskip

  A second purpose of the paper is to prove the following isometry extension result which 
is at least conceptually closely related to Theorem 1.2. However, while Theorem 1.2 is 
valid locally, this result depends crucially on global properties. 

\begin{theorem} \label{t1.3}
Let $g$ be a $C^{2}$ conformally compact Einstein metric on a compact manifold $M$ with 
$C^{\infty}$ boundary metric $(\partial M, \gamma)$, and suppose
\begin{equation}\label{e1.5}
\pi_{1}(M, \partial M) = 0.
\end{equation}
Then any connected group of isometries of $(\partial M, \gamma)$ extends to an action 
by isometries on $(M, g)$. 
\end{theorem}

  The condition (1.5) is equivalent to the statement that $\partial M$ is connected 
and  the inclusion map $\iota: \partial M \rightarrow M$ induces a surjection
$\pi_{1}(\partial M) \rightarrow \pi_{1}(M) \rightarrow 0$.

  Rather surprisingly, this result is closely related to the equations at conformal 
infinity induced by the Gauss-Codazzi equations on hypersurfaces tending to $\partial M$. 
It turns out that isometry extension from the boundary at least into a thickening of 
the boundary is equivalent to the requirement that the Gauss-Codazzi equations induced 
at $\partial M$ are preserved under arbitrary deformations of the boundary metric. This 
is discussed in detail in \S 5, see e.g.~Proposition 5.4. We note that this result does 
not hold for complete, asymptotically (locally) flat Einstein metrics, cf. Remark 5.8. 

  A simple consequence of Theorem 1.3 is the following uniqueness result: 
\begin{corollary} \label{c1.4}
A $C^{2}$ conformally compact Einstein metric with conformal infinity given by the 
class of the round metric $g_{+1}$ on the sphere $S^{n}$ is necessarily isometric 
to the Poincar\'e metric on the ball $B^{n+1}$. 
\end{corollary}

  Results similar to Theorem 1.3 and Corollary 1.4 have previously been proved in a number 
of different special cases by several authors, see for example [7], [9], [31], [33]; 
the proofs in all these cases are very different from the proof given here. 

\medskip

  It is well-known that unique continuation does not hold for large classes of elliptic 
systems of PDE's, even for general small perturbations of systems which are diagonal at 
leading order; see for instance [23] and references therein for a discussion related to 
geometric PDEs. The proofs of Theorems 1.1 and 1.2 rely on unique continuation results 
of Calder\'on [13], [14] and Mazzeo [27] respectively, based on Carleman estimates. The 
main difficulty in reducing the proofs to these results is the diffeomorphism covariance 
of the Einstein equations and, more importantly, that of the ``abstract'' Cauchy data 
$(\gamma, A)$ or $(\gamma, g_{(n)})$ at $\partial M$. The unique continuation theorem 
of Mazzeo requires a diagonal (i.e.~uncoupled) Laplace-type system of equations, at 
leading (second) order. The unique continuation result of Calder\'on is more general, 
but again requires strong restrictions on the structure of the leading order symbol 
of the operator. For emphasis and clarity, these issues are discussed in more detail 
in \S 2. The proofs of Theorems 1.1, 1.2 and 1.3 are then given in \S 3, \S 4 and \S 5 
respectively. 

\medskip

   I would like to thank the referee for very constructive remarks and comments on 
an earlier version of this paper, as well as Michael Taylor for interesting discussions 
on geodesic-harmonic coordinates and P. Chru\'sciel and E. Delay for interesting 
discussions concerning Theorem 1.3.

\section{Local Coordinates and Cauchy Data}
\setcounter{equation}{0}

 In this section, we discuss in more detail the remarks in the Introduction 
on classes of local coordinate systems, and their relation with Cauchy data 
on the boundary $\partial M$. 

  Thus, consider for example solutions to the system
\begin{equation} \label{e2.1}
Ric_{g} = 0, 
\end{equation}
defined near the boundary $\partial M$ of an $(n+1)$-dimensional manifold $M$. 
Since the Ricci curvature involves two derivatives of the metric, 
Cauchy data at $\partial M$ consist of the boundary metric $\gamma $ 
and its first derivative, invariantly represented by the $2^{\rm nd}$ 
fundamental form $A$ of $\partial M$ in $M$. Thus, we assume $(\gamma , A)$ 
are prescribed at $\partial M$, (subject to the Gauss and Gauss-Codazzi equations), 
and call $(\gamma, A)$ abstract Cauchy data. Observe that the abstract Cauchy 
data are invariant under diffeomorphisms of $M$ equal to the identity at 
$\partial M$. 

\medskip

 The metric $g$ determines the geodesic defining function 
$$t(x) = dist_{g}(x, \partial M).$$
The function $t$ depends of course on $g$; however, given any 
other smooth metric $g'$, there is a diffeomorphism $F$ of a 
neighborhood of $\partial M$, equal to the identity on $\partial M$, 
such that $t'(x) = dist_{F^{*}g'}(x, \partial M)$ satisfies $t' = t$. 
As noted above, this normalization does not change the abstract Cauchy data 
$(\gamma, A)$ and preserves the isometry class of the metric. 

  Let $\{y^{\alpha}\}$, $0 \leq \alpha \leq n$, be any local coordinates 
on a domain $\Omega $ in $M$ containing a domain $U$ in $\partial M$. 
We assume that $\{y^{i}\}$ for $1 \leq i \leq n$ form local coordinates 
for $\partial M$ when $y^{0} = 0$, so that $\partial / \partial y^{0}$ is 
transverse to $\partial M$. Throughout the paper, Greek indices $\alpha$, 
$\beta$ run from $0$ to $n$, while Latin indices $i$, $j$ run from $1$ 
to $n$. 

  If $g_{\alpha\beta}$ are the components of $g$ in these coordinates, then 
the abstract Cauchy problem associated to (2.1) in the local coordinates 
$\{y^{\alpha}\}$ is the system 
\begin{equation} \label{e2.2}
(Ric_{g})_{\alpha\beta} = 0, \ \ {\rm with} \ \  
g_{ij}|_{U} = \gamma_{ij}, \ \tfrac{1}{2}({\mathcal L}_{\nabla t}g)_{ij}|_{U} = 
a_{ij},
\end{equation}
where $\gamma_{ij}$ and $a_{ij}$ are given on $U$, (subject to the constraints 
of the Gauss and Gauss-Codazzi equations). Here one immediately sees a 
problem, in that (2.2) on $U \subset \partial M$ involves only the tangential 
part $g_{ij}$ of the metric (at 0 order), and not the full metric $g_{\alpha\beta}$ 
at $U$. The normal $g_{00}$ and mixed $g_{0i}$ components of the metric are not 
prescribed at $U$. As seen below, these components are gauge-dependent; 
they cannot be prescribed ``abstractly'', independent of coordinates, as is the 
case with $\gamma$ and $A$.

  In other words, if (2.1) is expressed in local coordinates $\{y^{\alpha}\}$ 
as above, then a well-defined Cauchy or unique continuation problem has the form 
\begin{equation} \label{e2.3}
(Ric_{g})_{\alpha\beta} = 0, \ \ {\rm with} \ \  
g_{\alpha\beta} = \gamma_{\alpha\beta}, \ 
\tfrac{1}{2}\partial_{t}g_{\alpha\beta} = a_{\alpha\beta}, \ {\rm on} \  U \subset 
\partial M , 
\end{equation}
where $\Omega$ is an open set in $({\mathbb R}^{n+1})^{+}$ with $\partial \Omega = U$ an 
open set in $\partial ({\mathbb R}^{n+1})^{+} ={\mathbb R}^{n}$. Formally, (2.3) is 
a determined system, while (2.2) is underdetermined. 

\medskip

 Let $g_{0}$ and $g_{1}$ be two solutions to (2.1), with the same 
Cauchy data $(\gamma , A)$, and with geodesic defining functions 
$t_{0}$, $t_{1}$. Changing the metric $g_{1}$ by a diffeomorphism 
if necessary, one may assume that $t_{0} = t_{1}$. One may then write the metrics 
with respect to a Gaussian or geodesic boundary coordinate system $(t, y^{i})$ as 
\begin{equation} \label{e2.4}
g_{k} = dt^{2} + (g_{k})_{t}, 
\end{equation}
where $(g_{k})_{t}$ is a curve of metrics on $\partial M$ and $k = 0, 1$. 
Here $y_{i}$ are coordinates on $\partial M$ which are extended into $M$ to be 
invariant under the flow of the vector field $\nabla t$. The metric $(g_{k})_{t}$ 
is the metric induced on $S(t)$ and pulled back to $\partial M$ by the flow of 
$\nabla t$. One has $(g_{k})_{0} = \gamma$ and $\frac{1}{2}\frac{d}{dt}(g_{k})_{t}|_{t=0} 
= A$. 

  Since $g_{0\alpha} = \delta_{0\alpha}$ in these coordinates, $\nabla t = \partial_{t}$, 
and hence the local coordinates are the same for both metrics, (or at least may be 
chosen to be the same). Thus, geodesic boundary coordinates are natural from the point 
of view of the Cauchy or unique continuation problem, since in such local coordinates 
the system (2.2), together with the prescription $g_{0\alpha} = \delta_{0\alpha}$, 
is equivalent to the system (2.3).

  However, the Ricci curvature is not elliptic or diagonal to leading order 
in these coordinates. The expression of the Ricci curvature in such coordinates 
does not satisfy the hypotheses of Calder\'on's theorem [14], and it appears to 
be difficult to establish unique continuation of solutions in these coordinates.

  Next suppose that $\{x^{\alpha}\}$ are boundary harmonic coordinates, 
defined as follows. For $1 \leq  i \leq n$, let $\hat x^{i}$ be local harmonic 
coordinates on a domain $U$ in $(\partial M, \gamma)$. Extend $\hat x^{i}$ into 
$M$ to be harmonic functions in $(\Omega, g)$, $\Omega \subset  M$, with 
Dirichlet boundary data; thus
\begin{equation} \label{e2.5}
\Delta_{g}x^{i} = 0, \ \ x^{i}|_{U} = \hat x^{i}.
\end{equation}
Let $x^{0}$ be a harmonic function on $\Omega$ with 0 boundary data, so that
\begin{equation} \label{e2.6}
\Delta_{g}x^{0} = 0, \ \ x^{0}|_{U} = 0.
\end{equation}
Then the collection $\{x^{\alpha}\}$, $0 \leq \alpha \leq  n$, form a 
local harmonic coordinate chart on a domain $\Omega \subset (M, g)$. In 
such coordinates, one has
\begin{equation} \label{e2.7}
(Ric_{g})_{\alpha\beta} = -\tfrac{1}{2}g^{\mu\nu}\partial_{\mu}\partial_{\nu}g_{\alpha\beta} 
+ Q_{\alpha\beta}(g, \partial g), 
\end{equation}
where $Q(g, \partial g)$ depends only on $g$ and its first derivatives. This is an 
elliptic operator, diagonal at leading order, and satisfies the hypotheses of 
Calder\'on's theorem. However, in general, the local Cauchy problem (2.3) is not 
well-defined in these coordinates; if $g_{0}$ and $g_{1}$ are two solutions of (2.1), 
each with corresponding local boundary harmonic coordinates, then the components 
$(g_{0})_{0\alpha}$ and $(g_{1})_{0\alpha}$ in general will differ at $U \subset 
\partial M$. This is of course closely related to the fact that there are many possible 
choices of harmonic functions $x^{\alpha}$ satisfying (2.5) and (2.6), and to the fact 
that the behavior of harmonic functions depends on global properties of $(\Omega, g)$. 
In any case, it is not known how to set up a well-defined Cauchy problem in these 
coordinates for which one can apply standard unique continuation results. 

  Consider then geodesic-harmonic coordinates ``intermediate'' between geodesic 
boundary and boundary harmonic coordinates. Thus, let $t$ be the geodesic distance 
to $\partial M$ as above. Choose local harmonic coordinates $\hat x^{i}$ on 
$\partial M$ as before and extend them into $M$ to be harmonic on the level sets 
$S(t)$ of $t$, i.e.~locally on $S(t)$, 
\begin{equation} \label{e2.8}
\Delta_{U(t)}x^{i} = 0, \ \ x^{i}|_{\partial U(t)} = \hat x^{i}|_{\partial U(t)}; 
\end{equation}
here the boundary value $\hat x^{i}$ is the extension of $\hat x^{i}$ on $U$ 
into $M$ which is invariant under the flow $\phi_{t}$ of $\nabla t$, 
and $U(t) = \phi_{t}(U) \subset S(t)$. The functions $(t, x^{i})$ form a coordinate 
system in a neighborhood $\Omega$ in $M$ with $\Omega\cap\partial M = U$. 

  It is not difficult to prove that geodesic-harmonic coordinates preserve the 
Cauchy data, in the sense that the data (2.2) in such coordinates imply the data 
(2.3). However, the Ricci curvature is not an elliptic operator in the metric 
in these coordinates, nor is it diagonal at leading order; the main reason is 
that the mean curvature of the level sets $S(t)$ is not apriori controlled. 
So again, it remains an open question whether unique continuation can be proved 
in these coordinates. 

  Having listed these attempts which appear to fail, a natural choice of 
coordinates which do satisfy the necessary requirements are $H$-harmonic 
coordinates $(\tau , x^{i})$, whose $\tau$-level surfaces $\Sigma_{\tau}$ 
are of prescribed mean curvature $H$ and with $x^{i}$ harmonic on $\Sigma_{\tau}$. 
These coordinates were introduced by Andersson-Moncrief [8] to prove a 
well-posedness result for the Cauchy problem for the Einstein equations 
in general relativity, and, as shown in [8], have a number of advantageous 
properties. Thus, adapting some of the arguments of [8], we show in \S 3 that 
the Einstein equations (1.2) are effectively elliptic in such coordinates, and 
such coordinates preserve the Cauchy data in the sense above, (i.e.~(2.2) 
implies (2.3)). It will then be shown that unique continuation holds in 
such coordinates, via application of the Calder\'on theorem.

\section{Proof of Theorem 1.1}
\setcounter{equation}{0}

 Theorem 1.1 follows from a purely local result, which we formulate as 
follows. Let $\Omega$ be a domain diffeomorphic to a cylinder $I\times 
B^{n} \subset {\mathbb R}^{n+1}$, with $U = \partial\Omega$ diffeomorphic to 
a ball $B^{n} \subset {\mathbb R}^{n} = \{0\}\times {\mathbb R}^{n} \subset 
{\mathbb R}^{n+1}$. Let $g$ be a Riemannian metric on $\Omega$ which is 
$C^{2,\alpha}$ up to $\partial\Omega$ in a given coordinate system 
$\{y^{\alpha}\}$ with $y^{0} = 0$ on $\partial\Omega$. Without loss 
of generality, we assume that $\Omega$ is chosen sufficiently small so 
that $g$ is close to the Euclidean metric $\delta$ in the $C^{2,\alpha}$ 
topology. For simplicity, rescale $(\Omega, g)$ and the coordinates 
$\{y^{\alpha}\}$ if necessary so that $(\Omega, g)$ is $C^{2,\alpha}$ 
close to the standard cylinder $((I\times B^{n}(1), B^{n}(1))) \subset 
({\mathbb R}^{n+1},{\mathbb R}^{n})$, $I = [0,1]$. The full boundary of $\Omega$, 
i.e.~$\bar \Omega \setminus \Omega$ will be denoted by $\bar \partial \Omega$.

 We will prove the following local version of Theorem 1.1. 
\begin{theorem} \label{t 3.1.}
  Let $g_{0}$, $g_{1}$ be two $C^{2,\alpha}$ metrics as above on $\Omega$ 
satisfying
\begin{equation} \label{e3.1}
Ric_{g_{k}} = \lambda g_{k}, 
\end{equation}
for some fixed constant $\lambda$, $k = 0,1$. Suppose $g_{0}$ and 
$g_{1}$ have the same abstract Cauchy data on $U$ in the sense of \S 2, 
so that $\gamma_{0} = \gamma_{1}$ and $A_{0} = A_{1}$. 

 Then $(\Omega, g_{0})$ is isometric to $(\Omega, g_{1})$, by an isometry 
equal to the identity on $\partial \Omega$. In particular, Theorem 1.1 holds. 
\end{theorem}

 The proof of Theorem 3.1 will proceed in several steps, organized 
around several Lemmas. We first work with a fixed metric $g$ on $\Omega$ 
as above. Let $N$ be the inward unit normal to $\partial\Omega$ in 
$\Omega$ and let $A = \nabla N$ be the corresponding $2^{\rm nd}$ fundamental 
form, with mean curvature $H = tr A$ on $U$. By the initial assumptions on 
$\Omega$ above, one has $A \sim 0$ and $H \sim 0$ in $C^{1,\alpha}$. 

 To begin, we construct the system of $H$-harmonic coordinates discussed 
at the end of \S 2. 

\begin{lemma} \label{l 3.2.}
  The domain $\Omega$ has a $C^{2,\alpha}$ foliation by surfaces $\Sigma_{s}$ 
with fixed mean curvature $H$, given by the mean curvature of $\partial\Omega$. 
Thus,
\begin{equation} \label{e3.2}
H_{\Sigma_{s}} = H, 
\end{equation}
in the coordinates $y^{\alpha}$, where $\Sigma_{s}$ is a graph over $U 
= \partial\Omega $ of the form $\Sigma_{s} = \{(f_{s}(y^{i}), y^{i}\}$ 
such that $y^{0}(\Sigma_{s}\cap \bar \partial \Omega) = s$, $s \in [0,1]$. 
Hence, in the coordinates $y^{i}$ for $\Sigma_{s}$, $H_{\Sigma_{s}}(y^{i}) 
= H(y^{i})$. 
\end{lemma}

{\bf Proof:} For any given $s > 0$, one may solve the prescribed mean 
curvature equation for graphs over $U = \partial\Omega$ with respect 
to the coordinates $y^{\alpha}$ in the metric $g$. The solution 
$f_{s}(y^{i})$ of (3.2) with $f_{s}|_{\partial U} = s$ exists and is 
unique, since the prescribed data $H \in C^{1,\alpha}$ is small and 
the metric $g$ is $C^{2,\alpha}$ close to the Euclidean metric on the 
domain $\Omega$. Note here that $U$ is close to the unit ball 
$B^{n}\subset{\mathbb R}^{n}$ with boundary mean curvature close to $n-1$. 
Since $g \in  C^{2,\alpha}$ and $H \in C^{1,\alpha}$, the solution $f_{s} 
\in C^{3,\alpha}(y^{i})$. Of course $f_{0} = 0$, so that $\Sigma_{0} = 
\partial \Omega$. The existence result above follows from, 
for instance, [19, Thm.~16.10, Thm.~15.10], cf.~also [10], while the 
uniqueness statement follows from [19, Thm.~10.7(iii)]. 

 The uniqueness implies that the family $\{\Sigma_{s}\}$ forms a $C^{0}$ 
foliation of $\Omega$. The derivative $\frac{df_{s}}{ds}$ satisfies a 
linear elliptic equation, with boundary value 1 and elliptic regularity 
implies the foliation is $C^{2,\alpha}$ smooth. A basis for the tangent 
space of $\Sigma_{s}$ is given by the vectors $\partial_{y_{i}} + 
(\partial_{y_{i}}f)\partial_{y^{0}}$ and hence the induced metric 
$g_{ij}$ on each $\Sigma_{s}$ in the induced $y$-coordinates is 
$C^{2,\alpha}$. Further, the $2^{\rm nd}$ fundamental form $A_{s}$ 
of $\Sigma_{s}$ is $C^{1,\alpha}$ in the $y$-coordinates. 

{\endproof}

 Given the foliation $\Sigma_{s}$ of prescribed mean curvature, we 
construct a corresponding $H$-harmonic coordinate system essentially as 
in \S 2. Thus, define the coordinate function $\tau$ by $\tau^{-1}(s) 
= \Sigma_{s}$, so that the level sets of $\tau$ are the leaves 
$\Sigma_{s}$, with $\tau  = s = y^{0}$ on $\partial\Sigma_{s}$. In 
particular, $\partial\Omega  = \tau^{-1}(0)$. For the spatial or 
tangential coordinates $x^{i}$, as before, choose harmonic coordinates 
$x^{i}$ on $U = \partial\Omega$, $\Delta_{\gamma}x^{i} = 0$, such that 
in these coordinates $\gamma_{ij}$ is $C^{2,\alpha}$ close to 
$\delta_{ij}$ on $U$. Extend $x^{i}$ to coordinates on each 
$\Sigma_{\tau}$ by solving the Dirichlet problem
\begin{equation} \label{e3.3}
\Delta_{\Sigma_{\tau}}x^{i} = 0, \ \  x^{i}|_{\partial\Sigma_{\tau}} = 
\bar x^{i}|_{\partial\Sigma_{\tau}}; 
\end{equation}
here the boundary value $\bar x^{i}$ is the 'vertical' extension of 
$x^{i}$ on $U$ into $\Omega$, i.e.~$\bar x^{i} = \tau  = y^{0}$ on 
$\partial\Sigma_{\tau}$. We may assume, without loss of generality, 
that $\{x^{\alpha}\} = \{(\tau , x^{i})\}$ form a coordinate system for 
$\Omega$. 

 As is well-known, the metric $g_{ij}$ induced on each $\Sigma_{\tau}$ 
has optimal regularity properties in harmonic coordinates, cf.~[11] for 
instance. Thus, the metric $g_{ij}^{x}$ on $\Sigma_{\tau}$ in the $x$-coordinate 
chart, $g_{ij}^{x} = g(\partial_{x^{i}}, \partial_{x^{j}})$, is as smooth as 
the induced metric $g_{ij}^{y}$ on $\Sigma_{\tau}$ in the $y$-coordinate 
chart. Since, by the above, $g_{ij}^{y} \in  C^{2,\alpha}$, one has 
$g_{ij}^{x}\in C^{2,\alpha}$, and $x^{i}\in C^{3,\alpha}(\{y^{j}\})$ 
on each $\Sigma_{\tau}$. 

 In the $\{y^{\alpha}\}$ coordinates, $\partial_{y^{0}}H_{\Sigma_{\tau}} = 0$ 
and $\partial_{y^{i}}H_{\Sigma_{\tau}}$ is fixed, independent of $\tau$. 
In the $\{x^{i}\}$ coordinates, one thus has on $\Sigma_{\tau}$,
\begin{equation} \label{e3.4}
\partial_{x^{i}}H = \frac{\partial y^{j}}{\partial x^{i}}\partial_{y^{j}}H  \ \ 
{\rm and} \ \ \partial_{\tau}H = \frac{\partial y^{i}}{\partial \tau}\partial_{y^{i}}H 
\in  C^{\alpha}(\{x^{k}\}). 
\end{equation}

 The metric $g = g_{\alpha\beta}$ in the $H$-harmonic coordinates 
$(\tau , x^{i})$ has the form
\begin{equation} \label{e3.5}
g = u^{2}d\tau^{2} + g_{ij}(dx_{i}+\sigma^{i}d\tau)(dx_{j}+\sigma^{j}d\tau), 
\end{equation}
where $\sigma  = \sigma^{i}\partial_{x^{i}}$ is the shift vector and 
$u$ is the lapse function, cf.~also [8]. One has $g^{00} = u^{-2}$, and $g_{00} = 
u^{2} + |\sigma|^{2}$, where $|\sigma|^{2} = g_{ij}\sigma^{i}\sigma^{j}$. 
The unit normal $N$ to the $\Sigma_{\tau}$ foliation is given by $N = 
\nabla\tau /|\nabla\tau|, \nabla\tau  = g^{0\alpha}\partial_{\alpha}$, 
so that $|\nabla\tau| = u^{-1}$ and 
\begin{equation} \label{e3.6}
N = u^{-1}(\partial_{\tau} +\sigma ), 
\end{equation}
with $\sigma  = u^{2}g^{0i}\partial_{i}$. In particular, $N(x^{i}) = 
ug^{0i}$, or equivalently,
\begin{equation} \label{e3.7}
u^{2}\nabla\tau (x^{i}) = \sigma^{i}. 
\end{equation}
In these coordinates, the $2^{\rm nd}$ fundamental form $A = 
\frac{1}{2}{\mathcal L}_{N}g$ of the leaves $\Sigma_{\tau}$ has the form
\begin{equation}\label{e3.8}
A_{ij} = {\tfrac{1}{2}}u^{-1}(\partial_{\tau}g_{ij} + 
({\mathcal L}_{\sigma}g)_{ij}),
\end{equation}
so that
\begin{equation} \label{e3.9}
\partial_{\tau}g_{ij} = 2uA_{ij} - ({\mathcal L}_{\sigma}g)_{ij}. 
\end{equation}

 A standard computation from commuting derivatives gives the Riccati equation 
\begin{equation} \label{e3.10}
({\mathcal L}_{N}A)_{ij} = A^{2}_{ij} - u^{-1}(D^{2}u)_{ij} - R_{N}(ij), 
\end{equation}
where $R_{N}(ij) = \langle R(\partial_{i},N)N, \partial_{j} \rangle$ 
and $A^{2}_{ij} = A_{i}^{k}A_{j}^{k}$. (The equation (3.10) may also be derived 
from the $2^{\rm nd}$ variation formula). Using (3.6) and the fact that $A$ 
is tangential, (i.e.~$A(N, \cdot ) = 0$), this gives
\begin{equation} \label{e3.11}
\partial_{\tau}A = -{\mathcal L}_{\sigma}A - D^{2}u + uA^{2} - uR_{N}. 
\end{equation}
Another straightforward calculation via the Gauss equations shows that 
$R_{N} = Ric_{g} - Ric_{\Sigma_{\tau}} + HA - A^{2}$, which, via 
(3.9) and (3.11) gives the system of 'evolution' equations for $g = g_{ij}$ 
and $A = A_{ij}$:
\begin{equation} \label{e3.12}
\partial_{\tau}g = 2uA - {\mathcal L}_{\sigma}g, 
\end{equation}
\begin{equation} \label{e3.13}
\partial_{\tau}A = -{\mathcal L}_{\sigma}A - D^{2}u + u[Ric_{\Sigma_{\tau}} 
- Ric_{g} + 2A^{2} - HA]. 
\end{equation}
(Up to sign differences, these are the well-known Einstein evolution 
equations in general relativity, cf.~[8], [32]). Substituting (3.12) in (3.13) 
gives the $2^{\rm nd}$ order evolution equation for $g$:
\begin{equation}\label{e3.14}
\partial_{\tau}^{2}g_{ij} = 2u\partial_{\tau}A_{ij} - 
({\mathcal L}_{\partial_{\tau}}{\mathcal L}_{\sigma}g)_{ij} + 
2(\partial_{\tau}u)A_{ij} 
\end{equation}
$$= -2u({\mathcal L}_{\sigma}A)_{ij} - 2u(D^{2}u)_{ij} + 
2u^{2}[Ric_{\Sigma_{\tau}} - Ric_{g} + 2A^{2} - HA]_{ij} - 
({\mathcal L}_{\partial_{\tau}}{\mathcal L}_{\sigma}g)_{ij} + 
2(\partial_{\tau}u)A_{ij}. $$
Since in $g$-harmonic coordinates $(Ric_{\Sigma_{\tau}})_{ij} = 
-\frac{1}{2}\Delta g_{ij} + Q_{ij}(g, \partial g)$, this implies, 
(via (3.8)),
\begin{equation} \label{e3.15}
(\partial_{\tau}^{2}+ u^{2}\Delta )g_{ij} = 
-2u^{2}(Ric_{g})_{ij}-2u(D^{2}u)_{ij} - 
2({\mathcal L}_{\sigma}{\mathcal L}_{\partial_{\tau}}g)_{ij} - 
({\mathcal L}_{\sigma}{\mathcal L}_{\sigma}g)_{ij} + Q_{ij}(g,\partial g), 
\end{equation}
where $Q_{ij}$ is a term involving at most the first order derivatives 
of $g_{\alpha\beta}$ in all $x^{\alpha}$ directions. Thus, for Einstein 
metrics (1.2),
\begin{equation} \label{e3.16}
(\partial_{\tau}^{2}+ u^{2}\Delta  + 
2\partial_{\sigma}\partial_{\tau} + 2\partial_{\sigma}^{2})g_{ij} = - 
2u(D^{2}u)_{ij} + Q_{ij}(g,\partial g), 
\end{equation}
where $Q_{ij}$ has the same general form as before.

 One also has the 'constraint' equations along each leaf 
$\Sigma_{\tau}$, involving the non-tangential part of the Ricci 
curvature:
\begin{equation} \label{e3.17}
\delta (A - Hg) = Ric(N, \cdot  ) = 0, 
\end{equation}
$$|A|^{2} - H^{2} + R_{\Sigma_{\tau}} = R_{g} - 2Ric_{g}(N,N) = 
(n-1)\lambda . $$

 Next, we derive the equations for the lapse $u$ and shift $\sigma$. 

\begin{lemma} \label{l 3.3.}
  The lapse $u$ and shift $\sigma$ satisfy the following equations:
\begin{equation} \label{e3.18}
\Delta u + |A|^{2}u + \lambda u = -uN(H) = 
-(\partial_{\tau}+\sigma)H. 
\end{equation}
\begin{equation} \label{e3.19}
\Delta\sigma^{i} = 2u\langle D^{2}x^{i}, A \rangle + u\langle dx^{i}, 
dH \rangle + 2\langle dx^{i}, A(\nabla u) - {\tfrac{1}{2}}Hdu \rangle . 
\end{equation}
\end{lemma}

{\bf Proof:} The lapse equation is derived by taking the trace of (3.10), 
and noting that $tr {\mathcal L}_{N}A = N(H) + 2|A|^{2}$. 

 For the shift equation, since the functions $x^{i}$ are harmonic on 
$\Sigma_{\tau}$, one has
$$\Delta ((x^{i})') + (\Delta')(x^{i}) = 0, $$
where $'$ denotes the Lie derivative with respect to $u^{2}\nabla\tau  
= uN$ and the Laplacian is taken with respect to the induced metric on 
the slices $\Sigma_{\tau}$. Thus, by (3.7), $(x^{i})' = \sigma^{i}$. 

 From standard formulas, cf.~[11, Ch.~1K] for example, one has
$$(\Delta')(x^{i}) = -2\langle D^{2}x^{i}, \delta^{*}uN\rangle 
+ 2\langle dx^{i}, \beta (\delta^{*}uN)\rangle, $$
where all the terms on the right are along $\Sigma_{\tau}$ and $\beta$ 
is the Bianchi operator, $\beta(k) = \delta k + \frac{1}{2}dtr k$. 
Thus, $\delta^{*}uN = uA$, and the shift components $\sigma^{i}$ satisfy
$$\Delta\sigma^{i} = 2u\langle D^{2}x^{i}, A\rangle - 2u\langle dx^{i}, 
\delta A + {\tfrac{1}{2}}dH\rangle + 2\langle dx^{i}, A(\nabla u) - 
{\tfrac{1}{2}}Hdu \rangle.$$
The relation (3.19) then follows from the constraint equation 
(3.17).

\begin{lemma} \label{l 3.4.}
  In local $H$-harmonic coordinates, the lapse-shift coordinates $(u,\sigma)$ 
are uniquely determined by the tangential metric $g_{ij}$ and $2^{\rm nd}$ 
fundamental form $A_{ij}$ on each $\Sigma_{\tau}$. Moreover, the abstract 
Cauchy data $(g_{ij}, A_{ij})$ on each leaf $\Sigma_{\tau}$ determine the local 
Cauchy data {\rm (2.2)}.
\end{lemma}

{\bf Proof:}  The system (3.18)-(3.19) is a coupled elliptic system in the pair 
$(u, \sigma)$ on $\Sigma_{\tau}$, with boundary values on $\partial\Sigma_{\tau}$ 
given by 
\begin{equation} \label{e3.20}
u|_{\partial\Sigma_{\tau}} = \tau ,  \ \ \sigma|_{\partial\Sigma_{\tau}} = 0. 
\end{equation}
Consider first the equations (3.18)-(3.19) in the $y^{i}$ coordinates. 
Elliptic regularity, (cf.~[19,~Ch.~6] for instance), gives bounds for $u$ 
and $\sigma $ in $C^{2,\alpha}$ provided the following terms are all bounded 
in $C^{\alpha}(y^{i})$: $g_{ij}^{y}$, $\partial_{y}g_{ij}^{y}$, $|A|^{2}$, 
$D^{2}x_{i}$, $\partial_{x_{i}}H$, $\partial_{\tau}H$, and the 
coefficients of $\sigma^{i}$ in $\sigma(H)$. By Lemma 3.2, $g_{ij}^{y} \in 
C^{2,\alpha}$ and $A \in C^{1,\alpha}$. Also, $C^{\alpha}$ bounds on the 
Hessian $D^{2}x^{i}$ follow from elliptic regularity for the harmonic 
functions $x^{i}$. The bounds on the derivatives of $H$ follow from (3.4). 
Thus, all the coefficients of (3.18)-(3.19) are bounded in $C^{\alpha}$. 

 Since the metric $g_{ij}$ is close to the flat metric, (in the 
$C^{2,\alpha}$ topology), it is standard that there is then a unique 
solution to the elliptic boundary value problem (3.18)-(3.19)-(3.20). 
The solution $(u, \sigma)$ is uniquely determined by the coefficients 
$(g_{ij}^{y}, A_{ij}^{y})$, and the terms or coefficients containing derivatives 
of $H$ and derivatives of $x^{i}$ in (3.18)-(3.19). We claim that these 
are also uniquely determined by $(g_{ij}^{y}, A_{ij}^{y})$.

  First, regarding the term $N(H)$, write $N = \omega \partial_{y_{0}} + 
\xi \partial_{y^{i}}$, so that $(\omega, \xi)$ are the lapse-shift of the 
foliation $\Sigma_{\tau}$ with respect to the $y^{\alpha}$-coordinates. As 
preceding (3.4), the derivatives $\partial_{y^{\alpha}}H$ are fixed, so 
that $N(H)$ depends only on $(\omega, \xi)$, (to zero order). Hence $N(H)$ 
can be expressed uniquely in terms of $(g_{ij}^{y}, u, \sigma)$, with 
$(u,\sigma)$ inserted into the equation (3.18) on the right side. In addition, 
the functions $x^{k} = x^{k}(y^{\ell})$ are uniquely determined by the metric 
$g_{ij}^{y}$, and so, for instance $\partial_{x_{i}}H$ is uniquely determined 
by $g_{ij}^{y}$ via (3.4).  

  Combining the facts above, it follows that $(u, \sigma)$ is uniquely 
determined by $(g_{ij}^{y}, A_{ij}^{y})$. The transformation $\{y^{\ell}\} 
\rightarrow \{x^{k}\}$ is invertible and $(g_{ij}^{y}, A_{ij}^{y})$ is thus 
uniquely determined by $(g_{ij}^{x}, A_{ij}^{x}) \equiv (g_{ij}, A_{ij})$. 
Clearly, $(u, \sigma)$ depends smoothly on $(g_{ij}, A_{ij})$. Of course, 
by (3.5) the data $(g_{ij}, u, \sigma)$ determine the full metric 
$g = g_{\alpha\beta}$. 

  Next we claim that $\partial_{\tau}g_{0\alpha}$ is also determined by 
$(g_{ij}, A_{ij})$ along $\Sigma_{\tau}$. To see this, one has 
$\frac{1}{2}({\mathcal L}_{N}g)_{0i} = A(\partial_{\tau}, \partial_{i}) = 
A(uN - \sigma, \partial_{i}) = -\sigma_{k}A_{ik}$, and similarly, 
$\frac{1}{2}({\mathcal L}_{N}g)_{00} = \sigma_{j}\sigma_{i}A_{ij}$. 
Hence $({\mathcal L}_{N}g)_{0\alpha}$ is determined by $(g_{ij}, A_{ij})$. 
Expanding this Lie derivative using (3.6), a simple computation gives
\begin{equation}\label{e3.21}
\partial_{\tau}g_{0i} - \partial_{\tau}(\log u)(g_{0i} + \sigma^{k}g_{ik}) = 
\phi_{1},
\end{equation}
$$\partial_{\tau}(\log |\sigma|) - \partial_{\tau}(\log u) = \phi_{2},$$
where $\phi_{1}$, $\phi_{2}$ are determined by $(g_{ij}, A_{ij})$. This system 
of equations is uniquely solvable for $\partial_{\tau}g_{0\alpha}$, showing 
that $\partial_{\tau}g_{0\alpha}$ is determined by $(g_{ij}, A_{ij})$. 
In particular, it follows that $(\gamma, A)$ determines the full Cauchy data 
$(g_{\alpha\beta}, \partial_{\tau}g_{\alpha\beta})$ at $U = \partial \Omega$. 

{\endproof}

  Summarizing the work above, the Einstein equations in local $H$-harmonic 
coordinates imply the following system on the data $(g_{ij}, u, \sigma)$:
\begin{equation} \label{e3.22}
(\partial_{\tau}^{2}+ u^{2}\Delta  + 
2\partial_{\sigma}\partial_{\tau} + 2u\partial_{\sigma}^{2})g_{ij} = - 
2u(D^{2}u)_{ij} + Q(g_{\alpha\beta}, \partial g_{\alpha\beta}), 
\end{equation}
\begin{equation} \label{e3.23}
\Delta u + |A|^{2}u + \lambda u = -(\partial_{\tau}+\sigma)H. 
\end{equation}
\begin{equation} \label{e3.24}
\Delta\sigma^{i} = 2u\langle D^{2}x^{i}, A \rangle  + u\partial_{x^{i}}H + 
2\langle dx^{i}, A(\nabla u) - {\tfrac{1}{2}}Hdu \rangle. 
\end{equation}

\begin{remark} \label{r 3.5.}
{\rm  The system (3.22)-(3.24) is essentially an elliptic system in 
$(g_{ij}, u, \sigma)$, given that $H$ is prescribed. Thus, assuming $u 
\sim 1$ and $\sigma \sim 0$, the operator $P = \partial_{\tau}^{2}+ 
u^{2}\Delta  + 2\partial_{\sigma}\partial_{\tau} + 
2u\partial_{\sigma}^{2}$ is elliptic on $\Omega$ and acts diagonally on 
$\{g_{ij}\}$, as is the Laplace operator on the slices $\Sigma_{\tau}$ 
acting on $(u, \sigma)$. The system (3.22)-(3.24) is of course coupled, 
but the couplings are all of lower order, i.e.~$1^{\rm st}$ order, except 
for the term $D^{2}u$ in (3.22). However, this term can be controlled or 
estimated by elliptic regularity applied to the lapse equation (3.23), 
(as discussed further below). 

 Given the above, it is not difficult to deduce that local $H$-harmonic 
coordinates have the optimal regularity property, i.e.~if $g$ is in 
$C^{m,\alpha}(\Omega)$ in some local coordinate system, then $g$ is in 
$C^{m,\alpha}(\Omega)$ in $H$-harmonic coordinates. Since this will not 
actually be used here, we omit further details of the proof. }
\end{remark}

{\bf  Proof of Theorem 3.1.}

\medskip

 Lemma 3.4 implies that it suffices to prove the unique continuation 
property for the tangential metric $g_{ij}$ in $H$-harmonic coordinates, 
since $g_{ij}$ determines the full metric $g_{\alpha\beta} = (g_{ij}, 
u, \sigma)$. 

 Thus suppose $g_{0}$ and $g_{1}$ are two Einstein metrics on $\Omega$ 
with identical $(\gamma, A)$ on $\partial\Omega$. One may construct 
$H$-harmonic coordinates for each of $g_{0}$ and $g_{1}$, and via a 
diffeomorphism identifying these coordinates, assume that the resulting 
pair of metrics $g$ and $\widetilde g$ have fixed $H$-harmonic coordinates 
$(\tau, x^{i})$, and both metrics satisfy the system (3.22)-(3.24). 
Let
\begin{equation} \label{e3.25}
h = h_{ij} = \widetilde g_{ij}- g_{ij}. 
\end{equation}
One then takes the difference of (3.22) and freezes the coefficents at 
$g$ to obtain a linear equation in $h$. Thus, for example, 
$\Delta_{\widetilde g}\widetilde g_{ij} - \Delta_{g}g_{ij} = 
\Delta_{g}(h_{ij}) - (g^{ab} - \widetilde g^{ab})
\partial_{a}\partial_{b}\widetilde g_{ij}$. The second 
term here is of zero order, (rational), in the difference $h$, with 
coefficients depending on two derivatives of $\widetilde g$. Carrying out 
the same procedure on the remaining terms in (3.22) gives the equation
$$(\partial_{\tau}^{2}+ u^{2}\Delta  + 
2\partial_{\sigma}\partial_{\tau} + 2u\partial_{\sigma}^{2})h_{ij} = - 
2(\widetilde u(\widetilde D^{2}\widetilde u)_{ij} - u(D^{2}u)_{ij}) + 
Q_{ij}(h_{\alpha\beta},\partial_{\mu} h_{\alpha\beta}),$$
where $Q$ depends on two derivatives of the background $\widetilde g$, but 
only on one derivative, (in all directions), of $h_{\alpha\beta}$. Similarly, 
$\widetilde D^{2}\widetilde u$ - $D^{2}u = D^{2}v + (\widetilde D^{2} - 
D^{2})\widetilde u$, where $v = \widetilde u - u$ and the second term is 
of the form $Q$ above. Hence, 
\begin{equation} \label{e3.26}
(\partial_{\tau}^{2}+ u^{2}\Delta  + 
2\partial_{\sigma}\partial_{\tau} + 2u\partial_{\sigma}^{2})h_{ij} = - 
2u(D^{2}v)_{ij} + Q_{ij}(h_{\alpha\beta},\partial_{\mu} h_{\alpha\beta}), 
\end{equation}
Note that since we have linearized, $Q$ depends linearly on $h_{\alpha\beta}$ 
and $\partial_{\mu} h_{\alpha\beta}$, with nonlinear coefficients depending 
on $\widetilde g$ and $g$. 

 Next we use the lapse and shift equations (3.23)-(3.24) to estimate the 
differences $v = \widetilde u - u$ and $\chi = \widetilde \sigma - \sigma$. 
Thus, as before, $\Delta_{\widetilde g}\widetilde u - \Delta_{g}u = 
\Delta_{g}v + D^{2}_{h}(\widetilde u)$, where $D^{2}_{h}$ is a $2^{\rm nd}$ 
order differential operator on $\widetilde u$ with coefficients depending on 
the difference $h$, to zero order. Similarly, $|\widetilde A|^{2}\widetilde u - 
|A|^{2}u = |A|^{2}v + D^{0}_{h}(\widetilde u)$, the latter depending on $h$ to 
first order. As in the proof of Lemma 3.4, the term $N(H)$ depends only on 
$(g_{ij}, u, \sigma)$ to zero order. Taking the difference, it then follows 
from (3.23) and Lemma 3.4 that 
\begin{equation} \label{e3.27}
\Delta v + |A|^{2}v + \lambda v = Q(h_{ij}, \partial_{\mu}h_{ij}, v, \chi), 
\end{equation}
with $Q$ depending on two derivatives of $\widetilde g$ and linear in its 
arguments. Since $v = 0$ on $\partial\Sigma_{\tau}$, elliptic regularity then 
gives
\begin{equation} \label{e3.28}
||v||_{L_{x}^{2,2}} \leq  C(\widetilde g, g)[||h_{ij}||_{L_{(\tau ,x)}^{1,2}} 
+ ||\chi||_{L_{x}^{2}}], 
\end{equation}
where the $L^{2,2}$ norm on the left is with respect to spatial 
derivatives $x$, while the $L^{1,2}$ norm on the right includes also the 
time derivative $\tau$; both norms are taken along the leaves $\Sigma_{\tau}$. 
Working in the same way with the shift equation (3.24) gives the analogous 
estimate
\begin{equation} \label{e3.29}
||\chi||_{L_{x}^{2,2}} \leq  C(\widetilde g, g)[||h_{ij}||_{L_{(\tau ,x)}^{1,2}} 
+ ||v||_{L_{x}^{2}}]. 
\end{equation}

  It follows from (3.26) and (3.28) that
\begin{equation} \label{e3.30}
||P(h_{ij})||_{L_{x}^{2}} \leq  C(\widetilde g, g)
||h_{\alpha\beta}||_{L_{(\tau ,x)}^{1,2}}, 
\end{equation}
where $P$ is given as following (3.24). We claim that the right side of (3.30) 
is bounded by $||h_{ij}||_{L_{(\tau ,x)}^{1,2}}$. To see this, $h_{\alpha\beta} 
= (h_{ij}, v, \chi)$, so that it suffices to show that $v$, $\partial_{\tau}v$ 
and $\chi$, $\partial_{\tau}\chi$ are determined by $(h_{ij}, \partial_{\tau}h_{ij})$, 
with corresponding bounds on the $L^{2}$ norms. This follows directly from the proof 
of Lemma 3.4. Thus, $u$ and $\sigma$ are uniquely determined by $(g_{ij}, A_{ij})$, 
and hence by $(g_{ij}, \partial_{\tau}g_{ij})$, which implies that $v$ and $\chi$ 
are determined by linear equations in $(h_{ij}, \partial_{\tau}h_{ij})$; this 
also follows from (3.28)-(3.29). Similarly, since $\partial_{\tau}g_{0\alpha}$ 
is determined by $(g_{ij}, A_{ij})$ as in (3.21), it follows that 
$\partial_{\tau}h_{0\alpha}$ is determined by $(h_{ij}, \partial_{\tau}h_{ij})$. 
Hence, (3.30) gives 
\begin{equation} \label{e3.31}
||P(h_{ij})||_{L_{x}^{2}} \leq  C(\widetilde g, g)
||h_{ij}||_{L_{(\tau ,x)}^{1,2}}, 
\end{equation}

   As noted above, the operator $P$ is elliptic and diagonal, and by 
construction, the Cauchy data for $P$ vanish at $U = \partial\Omega$, 
i.e.
\begin{equation} \label{e3.32}
h = \partial_{\tau}h = 0 \ \ {\rm at} \ \  U, 
\end{equation}
since the lapse-shift $(u, \sigma)$ of $g$ and $\widetilde g$ and their 
$\tau$-derivatives agree at $U$. From now on, $h = h_{ij}$. 

 Now we claim that $P$ satisfies the hypotheses of the Calderon unique 
continuation theorem [14]. Following [14], decompose the symbol of $P$ 
as 
\begin{equation} \label{e3.33}
A_{2}(\tau ,x,\xi) = (u^{2}g^{kl}\xi_{k}\xi_{l} + 
2u\sigma^{k}\sigma^{l}\xi_{k}\xi_{l})I, 
\end{equation}
$$A_{1}(\tau ,x,\xi) = 2\sigma^{k}\xi_{k}I, $$
where $I$ is the $N\times N$ identity matrix, $N = \frac{1}{2}n(n+1)$, 
equal to the cardinality of $\{ij\}$. Setting $|\xi|^{2} = 1$, (3.33) 
becomes
$$A_{2}(\tau ,x,\xi) = (u^{2}+ 2u\sigma^{k}\sigma^{l}\xi_{k}\xi_{l})I,$$
$$A_{1}(\tau ,x,\xi) = 2\sigma^{k}\xi_{k}I.$$
Now form the matrix
\begin{equation}\label{e3.34}
M = \left(
\begin{array}{cc}
0  &  -I \\
A_{2} & A_{1}
\end{array}
\right)
\end{equation}

 The matrices $A_{1}$ and $A_{2}$ are diagonal, and it is then easy to 
see that $M$ is diagonalizable, i.e.~has a basis of eigenvectors over 
${\mathbb C}$. This implies that $M$ satisfies the hypotheses of 
[14, Thm.~11(iii)], cf.~also [14, Thm.~4]. The bound (3.31) is substituted 
in the basic Carleman estimate of [14, Thm.~6], cf.~also [29, (6.1)], showing 
that $h_{ij}$ satisfies the unique continuation property. 

 It follows from (3.32) and the Calder\'on unique continuation theorem that 
$$h_{ij} = \widetilde g_{ij} - g_{ij} = 0, $$
in an open neighborhood $\Omega' \subset \Omega$. By Lemma 3.4, this 
implies 
$$\widetilde g_{\alpha\beta} = g_{\alpha\beta}, $$
in $\Omega'$, so that $g_{1}$ is isometric to $g_{0}$ in $\Omega'$. 
By construction, the isometry from $g_{0}$ to $g_{1}$ equals the identity 
on $\partial \Omega$. 

  This shows that the metric $g$ is uniquely determined in $\Omega'$, up to 
isometry, by the abstract Cauchy data on $\partial \Omega$. Since Einstein 
metrics are real-analytic in the interior in harmonic coordinates, a standard 
analytic continuation argument, (cf.~[25] for instance), then implies that 
$g$ is unique up to isometry everywhere in $\Omega$. In the context of 
Theorem 1.1, the same analytic continuation argument shows that the Cauchy 
data of $g$ on $U$ uniquely determines the topology of $M$ and $\partial M$, 
up to covering spaces, as well as the Cauchy data $(\gamma, A)$ on $\partial M$ 
outside $U = \partial \Omega$. This completes the proof of Theorems 3.1 and 1.1.

{\endproof}

  As an illustration, suppose $(M_{1}, g_{1})$ and $(M_{2}, g_{2})$ are a 
pair of Einstein metrics on compact manifolds-with-boundary and the Cauchy 
data for $g_{1}$ and $g_{2}$ agree on an open set $U$ of the boundary. 
Suppose $M_{i}$ are connected and the inclusion map of $\partial M_{i}$ 
into $M_{i}$ induces a surjection of fundamental groups, i.e.
\begin{equation}\label{e3.35}
\pi_{1}(\partial M_{i}) \rightarrow \pi_{1}(M_{i}) \rightarrow 0,
\end{equation}
for $i = 1,2$, so that every loop in $M_{i}$ is homotopic to a loop in 
$\partial M_{i}$. Then $M_{1}$ is diffeomorphic to $M_{2}$ and $g_{1}$ is 
isometric to $g_{2}$. 

\medskip

 We conclude this section with a discussion of generalizations of 
Theorem 1.1. First, one might consider the unique continuation problem 
for 
\begin{equation} \label{e3.36}
Ric_{g} = T, 
\end{equation}
where $T$ is a fixed symmetric bilinear form on $M$, at least $C^{\alpha}$ up 
to $\bar M$. However, this problem is not natural, in that is not covariant under 
changes by diffeomorphism. For metrics alone, the Einstein equation (1.2) is the 
only equation covariant under diffeomorphisms which involves at most 
the $2^{\rm nd}$ derivatives of the metric.  Nevertheless, the proof of Theorem 
1.1 shows that if $\widetilde g$ and $g$ are two solutions of (3.36) which have 
common $H$-harmonic coordinates near (a portion of) $\partial M$ on which 
$(\gamma, A) = (\widetilde \gamma, \widetilde A)$, then $\widetilde g$ is 
isometric to $g$ near $\partial M$.

   Instead, it is more natural to consider the Einstein equation coupled 
(covariantly) to other fields $\chi$ besides the metric; such 
equations arise naturally in many areas of physics. For example, $\chi$ 
may be a function on $M$, i.e.~a scalar field, or $\chi$ may be a 
connection 1-form (gauge field) on a bundle over $M$. We assume that the 
field(s) $\chi$ arise via a diffeomorphism-invariant Lagrangian ${\mathcal L} 
= {\mathcal L}(g, \chi)$, depending on $\chi$ and its first derivatives in 
local coordinates, and that $\chi$ satisfies field equations, i.e.~Euler-Lagrange 
equations, coupled to the metric. For example, for a free massive scalar field, 
the equation is the eigenfunction equation
\begin{equation} \label{e3.37}
\Delta_{g}\chi  = \mu\chi ,
\end{equation}
while for a connection 1-form, the equations are the Yang-Mills 
equations, (or Maxwell equations when the bundle is a $U(1)$ bundle):
\begin{equation} \label{e3.38}
dF = d^{*}F = 0,
\end{equation}
where $F$ is the curvature of the connection $\chi$. Associated to 
such fields is the stress-energy tensor $T = T_{\mu\nu}$; this is a 
symmetric bilinear form obtained by varying the Lagrangian for $\chi$ 
with respect to the metric, cf.~[22] for example. For the free massive 
scalar field $\chi$ above, one has
$$T = d\chi\cdot  d\chi  - \tfrac{1}{2}(|d\chi|^{2} + \mu\chi^{2})g, $$
while for a connection 1-form 
$$T = F \cdot F - \tfrac{1}{4}|F|^{2}g, $$
where $(F \cdot F)_{ab} = F_{ac}F_{bd}g^{cd}$. 

 When the part of the Lagrangian involving the metric to $2^{\rm nd}$ order 
only contains the scalar curvature, i.e.~the Einstein-Hilbert action, the 
resulting coupled Euler-Lagrange equations for the system $(g, \chi)$ 
are
\begin{equation} \label{e3.39}
Ric_{g} - \frac{R}{2}g = T, \ \ E_{g}(\chi) = 0. 
\end{equation}
By taking the trace, this can be rewritten as
\begin{equation} \label{e3.40}
Ric_{g} = \hat T = T - \frac{1}{n-1}tr_{g}T,  \ \ E_{g}(\chi) = 0. 
\end{equation}

 Here we assume $E_{g}(\chi)$ is a $2^{\rm nd}$ order elliptic system for 
$\chi$, with coefficients depending on $g$, as in (3.37) or (3.38), (the latter 
viewed as an equation for the connection). In case the field(s) $\chi$ have an 
internal symmetry group, as in the case of gauge fields, this will require 
a particular choice of gauge for $\chi$ in which the Euler-Lagrange equations 
become an elliptic system in $\chi$. It is also assumed that solutions $\chi$ of 
$E_{g}(\chi) = 0$ satisfy the unique continuation property; for instance $E_{g}$ satisfies 
the hypotheses of the Calder\'on theorem [14]. Theorem 1.1 now easily extends to cover 
(3.39) or (3.40).

\begin{proposition} \label{p3.3}
Let $M$ be a compact manifold with boundary $\partial M$. Then $C^{2,\alpha}$ solutions 
$(g, \chi)$ of (3.39) on $\bar M$ are uniquely determined, up to local isometry, by 
the Cauchy data $(\gamma, A)$ of $g$ and the Cauchy data $(\chi, \partial_{t}\chi)$ 
on an open set $U \subset \partial M$. 
\end{proposition}

{\bf Proof:}
 The proof is the same as the proof of Theorem 1.1. Briefly, via a suitable 
diffeomorphism equal to the identity on $\partial M$, one brings a pair of solutions 
of (3.39) with common Cauchy data into a fixed system of $H$-harmonic coordinates 
for each metric. As before, one then applies Calder\'on uniqueness to the resulting 
system (3.39) in the difference of the metrics and fields. Further details are left 
to the reader. 
{\endproof}

\section{Proof of Theorem 1.2.}
\setcounter{equation}{0}

 Let $g$ be a conformally compact metric on a compact $(n+1)$-manifold $M$ with boundary 
which has a $C^{2}$ geodesic 
compactification
\begin{equation} \label{e4.1}
\bar g = t^{2}g, 
\end{equation}
where $t(x) = dist_{\bar g}(x, \partial M)$. By the Gauss Lemma, one has the splitting
\begin{equation} \label{e4.2}
\bar g = dt^{2} + g_{t}, 
\end{equation}
near $\partial M,$ where $g_{t}$ is a curve of metrics on $\partial M$ 
with $g_{0} = \gamma$ the boundary metric. The curve $g_{t}$ is 
obtained by taking the induced metric the level sets $S(t)$ of $t$, 
and pulling back by the flow of $N = \bar \nabla t$. Note that if 
$r = -\log t$, then $g = dr^{2} + t^{-2}g_{t}$, so the integral curves 
of $\nabla r$ with respect to $g$ are also geodesics. Each choice of boundary 
metric $\gamma \in [\gamma]$ determines a unique geodesic defining 
function $t$.

 Now suppose $g$ is Einstein, so that (1.4) holds and suppose for the moment that 
$g$ is $C^{2}$ conformally compact with $C^{\infty}$ smooth boundary metric $\gamma$. 
Then the boundary regularity result of [16] implies that $\bar g$ is $C^{\infty}$ 
smooth when $n$ is odd, and is $C^{\infty}$ polyhomogeneous when $n$ is even. Hence, 
the curve $g_{t}$ has a Taylor-type series in $t$, called the Fefferman-Graham 
expansion [18]. The exact form of the expansion depends on whether $n$ is odd 
or even. If $n$ is odd, one has a power series expansion 
\begin{equation} \label{e4.3}
g_{t} \sim g_{(0)} + t^{2}g_{(2)} + \cdots + t^{n-1}g_{(n-1)} + t^{n}g_{(n)} 
+ \cdots , 
\end{equation}
while if $n$ is even, the series is polyhomogeneous,
\begin{equation} \label{e4.4}
g_{t} \sim g_{(0)} + t^{2}g_{(2)} + \cdots + t^{n}g_{(n)} + t^{n}\log t \ {\mathcal H} + 
\cdots . 
\end{equation}
In both cases, this expansion is even in powers of $t$, up to $t^{n}$. It is important 
to observe that the coefficients $g_{(2k)}$, $k \leq [n/2]$, as well as the coefficient 
${\mathcal H}$ when $n$ is even, are explicitly determined by the boundary metric $\gamma  = 
g_{(0)}$ and the Einstein condition (1.4), cf.~[18], [20]. For $n$ even, the series (4.4) 
has terms of the form $t^{n+k}(\log t)^{m}$. 

  For any $n$, the divergence and trace (with respect to $g_{(0)} = \gamma$) of $g_{(n)}$ 
are determined by the boundary metric $\gamma$; in fact there is a symmetric bilinear form 
$r_{(n)}$ and scalar function $a_{(n)}$, both depending only on $\gamma$ and its derivatives 
up to order $n$, such that 
\begin{equation} \label{e4.5}
\delta_{\gamma}(g_{(n)} + r_{(n)}) = 0, \ \ {\rm and} \ \ tr_{\gamma}(g_{(n)} + r_{(n)}) = 
a_{(n)}.
\end{equation}
For $n$ odd, $r_{(n)} = a_{(n)} = 0$. (The divergence-free tensor $g_{(n)} + r_{(n)}$ is 
closely related to the stress-energy of a conformal field theory on $(\partial M, \gamma)$, 
cf.~[17]). The relations (4.5) will be discussed further in \S 5. 

  However, beyond the relations (4.5), the term $g_{(n)}$ is not determined by $g_{(0)}$; 
it depends on the ``global'' structure of the metric $g$. The higher order coefficients 
$g_{(k)}$ of $t^{k}$ and coefficients $h_{(km)}$ of $t^{n+k}(\log t)^{m}$, are then 
determined by $g_{(0)}$ and $g_{(n)}$ via the Einstein equations. The equations (4.5) 
are constraint equations, and arise from the Gauss-Codazzi and Gauss and Riccati equations 
on the level sets $S(t) = \{x: t(x) = t\}$ in the limit $t \rightarrow 0$; this is also 
discussed further in \S 5. 

   In analogy to the situation in \S 3, the term $g_{(n)}$ corresponds to the $2^{\rm nd}$ 
fundamental form $A$ of the boundary, in that, modulo the constraints (4.5), it is freely 
specifiable as Cauchy data, and is the only such term depending on normal derivatives of 
the boundary metric.

\medskip

Suppose now $g_{0}$ and $g_{1}$ are two solutions of 
\begin{equation} \label{e4.6}
Ric_{g} + ng = 0, 
\end{equation}
with the same $C^{\infty}$ conformal infinity $[\gamma]$. Then there exist 
geodesic defining functions $t_{k}$ such that $\bar g_{k} = (t_{k})^{2}g_{k}$ 
have a common boundary metric $\gamma \in [\gamma]$, and both metrics are 
defined for $t_{k} \leq  \varepsilon$, for some $\varepsilon > 0$. 

  The hypotheses of Theorem 1.2, together with the discussion above concerning 
(4.3) and (4.4), then imply that 
\begin{equation} \label{e4.7}
|g_{1} - g_{0}| = o(e^{-nr}) = o(t^{n}), 
\end{equation}
where the norm is taken with respect to $g_{1}$, (or $g_{0}$). 

\medskip

  Given this background, we prove the following more general version of Theorem 1.2, 
analogous to Theorem 3.1. Let $\Omega$ be a domain diffeomorphic to $I\times 
B^{n}$, where $B^{n}$ is a ball in ${\mathbb R}^{n}$ with boundary $U = \partial 
\Omega$ diffeomorphic to a ball in ${\mathbb R}^{n} \simeq \{0\}\times 
{\mathbb R}^{n}$. 

\begin{theorem} \label{t4.1}
Let $g_{0}$ and $g_{1}$ be a pair of conformally compact Einstein metrics on 
a domain $\Omega$ as above. Suppose $g_{0}$ and $g_{1}$ have $C^{2,\alpha}$ 
geodesic compactifications, and (4.7) holds in $\Omega$. 

  Then $(\Omega, g_{0})$ is isometric to $(\Omega, g_{1})$, by an isometry equal 
to the identity on $\partial \Omega$. Hence, if $(M_{0}, g_{0})$ and $(M_{1}, g_{1})$ 
are conformally compact Einstein metrics on compact manifolds with boundary, and 
(4.7) holds on some open domain $\Omega$ in $M_{0}$ and $M_{1}$, then the manifolds 
$M_{0}$ and $M_{1}$ are diffeomorphic in some covering space of each and the lifted 
metrics $g_{0}$ and $g_{1}$ are isometric.
\end{theorem}

 The proof of Theorem 4.1 is very similar to that of Theorem 3.1. For 
clarity, we first prove the result in case the metrics $g_{i}$, $i = 
0,1$, have a common $C^{\infty}$ boundary metric $\gamma$ and then show 
how the proof can be extended to cover the more general case of metrics 
with less regularity. 

 By applying a diffeomorphism if necessary, one may assume that the 
metrics $g_{i}$ have a common geodesic defining function $t$ defined 
near $\partial\Omega$ and common geodesic boundary coordinates. By [16], 
the geodesically compactified metrics $\bar g_{i} = t^{2}g_{i}$ are 
$C^{\infty}$ and extend $C^{\infty}$ to $\partial\Omega$. It follows 
from the discussion of the Fefferman-Graham expansion following (4.5) 
that $g_{0}$ and $g_{1}$ agree to infinite order at $\partial U$, i.e.
\begin{equation} \label{e4.8}
k = g_{1} - g_{0} = O(t^{\nu}), 
\end{equation}
for any $\nu  < \infty$. Of course $k_{0\alpha} = 0$. 

 For the rest of the proof, we work in the setting of the compactified 
metrics $\bar g_{i}$. As in the proof of Theorem 3.1, we assume 
that the domain $\Omega$ is sufficiently small so that $(\Omega , 
\bar g_{i})$ is smoothly close to the flat metric on the standard 
cylinder $I\times B^{n}$, with $\bar A = 0$ on $U = \partial\Omega$. 
(Note that $g_{(1)} = 0$ in (4.3)-(4.4)). In particular, near 
$\partial \Omega$, $\bar H = O(t)$. One may construct a foliation 
$\Sigma_{\tau}$ with $\bar H_{\Sigma_{\tau}} = 0$ exactly as in Lemma 3.2, 
together with corresponding $H$-harmonic coordinates $(\tau , x^{i})$. 
All of the analysis carried out in \S 3 through to Lemma 3.4 carries 
over to this situation with only a single difference. Namely, for the 
term $Ric_{g}$ in (3.14) or (3.15), one now no longer has $Ric_{g} = 
\lambda g$, but instead the Ricci curvature $\bar Ric$ of the compactified 
metric $\bar g$. Using the facts that $Ric_{g} = -ng$ and the compactification 
$\bar g$ is geodesic, standard formulas for the behavior of Ricci curvature under 
conformal change give
\begin{equation} \label{e4.9}
\bar Ric = -(n-1)t^{-1}\bar D^{2}t  - t^{-1}\bar \Delta t \bar g. 
\end{equation}
One has $\bar D^{2}t = {\mathcal L}_{\nabla t}\bar g = O(t)$, 
(recall that $\bar g$ is $C^{\infty}$ up to $\partial\Omega$). If 
$(t, y^{i})$ are geodesic boundary coordinates, then $\partial_{x^{i}} 
= \sum(1-\varepsilon (\tau))\partial_{y^{j}} + \varepsilon (\tau)\nabla t$, 
where $\varepsilon(\tau) = O(\tau)$. Similarly, $\tau /t = 1 + 
\varepsilon(\tau)$. (The specific form of $\varepsilon (\tau)$ 
of course differs in each occurance above, but this is insignificant). 
Since $\bar D^{2}t$ vanishes on $\nabla t$, it follows from (4.9) 
that in the $x^{\alpha}$ coordinates on $\Sigma_{\tau}$,
\begin{equation} \label{e4.10}
\bar Ric_{ij} = -(n-1)(1-\varepsilon)^{2}t^{-1}({\mathcal L}_{\nabla t}
\bar g)_{ij}  - (1-\varepsilon)^{2}t^{-1}(\bar \Delta t \bar)g_{ij} + 
\varepsilon t^{-1}(\bar \Delta t)q_{ij}, 
\end{equation}
where $q_{ij}$ depends only on $\bar g_{0\alpha}$ to zero-order. Next 
$({\mathcal L}_{\nabla t}\bar g) = (1 - \varepsilon)\partial_{\tau}\bar g + 
\varepsilon (\tau)\partial_{x^{\alpha}}\bar g$ and similarly for the 
Laplace term in (4.10). Substituting (4.10) in (3.15), it follows that 
the analogue of (3.16) in this context is the 'evolution equation'
\begin{equation} \label{e4.11}
\tau^{2}(\partial_{\tau}^{2}+ u^{2}\Delta  + 
2\partial_{\sigma}\partial_{\tau} + 2\partial_{\sigma}^{2})g_{ij} = - 
2\tau^{2}u(D^{2}u)_{ij} + Q_{ij}(g,\tau\partial g), 
\end{equation}
where $Q_{ij}$ is a term involving $g_{\alpha\beta}$ with at most first order 
derivatives of the form $\tau\partial_{\alpha}$. Here and below, we 
drop the bar from the notation. 

 The lapse $u$ and shift $\sigma$ satisfy essentially the same 
equations as before, namely
\begin{equation} \label{e4.12}
\Delta u + |A|^{2}u - (t^{-1}\Delta t)u = 0, 
\end{equation}
\begin{equation} \label{e4.13}
\Delta\sigma^{i} = 2u\langle D^{2}x^{i}, A \rangle + + 2\langle dx^{i}, A(\nabla u) 
\rangle 
\end{equation}
Comparing with (3.18)-(3.19), one has here $H = 0$, with the $\lambda$ 
term in (3.17) replaced by $-t^{-1}\Delta t$. Lemma 3.4 holds 
as before, since $t^{-1}\Delta t$ is smooth up to $\partial \Omega$. 

 One now proceeds just as in the proof of Theorem 3.1, taking the 
difference of the equation (4.11) to obtain a linear equation on $h = 
\widetilde g - g$; (recall that the bars have been removed from the 
notation). Note that by (4.8), together with elliptic regularity applied 
to (4.12)-(4.13), as in the proof of Lemma 3.4, one has
\begin{equation} \label{e4.14}
h_{\alpha\beta} = O(t^{\nu}), 
\end{equation}
for all $\nu  < \infty$.  The estimates (3.28)-(3.31) hold as before.

 Let $P(h_{ij}) = \tau^{2}(\partial_{\tau}^{2}+ u^{2}\Delta  + 
2\partial_{\sigma}\partial_{\tau} + 2\partial_{\sigma}^{2})$. Then 
$P$ is a fully degenerate $2^{\rm nd}$ order elliptic operator, with 
smooth coefficients, and one has
$$||P(h_{ij})||_{L_{x}^{2}} \leq  C||h_{ij}||_{L_{\tau,x}^{1,2}},$$
where the $1^{\rm st}$ order derivatives on the right are of the form 
$\tau\partial$. Further, by (4.14), $h$ vanishes to infinite order at 
$\partial\Omega$. It then follows from a unique continuation theorem 
of Mazzeo, [27, Thm.~14], 
that
$$h_{ij} = 0 $$
in $\Omega' \subset \Omega$. The vanishing of $h = h_{\alpha\beta}$ in 
$\Omega$ then follows as before in the proof of Theorem 3.1. 

  Next suppose $g_{0}$ and $g_{1}$ have only a $C^{2,\alpha}$ geodesic 
compactification with a common boundary metric $\gamma$, but that 
(4.7) holds. All of the arguments above remain valid, except the 
infinite order vanishing property (4.8), and the corresponding (4.14), 
which are replaced by the statements $k = o(t^{n})$ and $h = o(t^{n})$ 
respectively. The unique continuation result in [27] per se, requires 
the infinite order decay (4.14). Thus, it suffices to show that (4.14) 
does in fact hold. 

 To do this, we first show that $k = O(t^{\nu})$ weakly, for all $\nu  
< \infty$. This will imply $h = O(t^{\nu})$ weakly, and the strong or 
pointwise decay (4.14) then follows from elliptic regularity.

  In geodesic boundary coordinates, the geodesic compactification of a 
conformally compact Einstein metric satisfies the equation
\begin{equation} \label{e4.15}
t\ddot g - (n-1)\dot g - 2Hg^{T}-2tRic_{S(t)} + tH\dot g - t(\dot g)^{2} = 0,
\end{equation}
where $\dot g$ is the Lie derivative of $g$ with respect to $\nabla t$, 
cf.~[18] or [21]. Thus $\dot g = 2A$, where $A$ is the $2^{\rm nd}$ 
fundamental form of the level set $S(t)$ of $t$, (with respect to the 
inward normal). Also $H = tr A$, $T$ denotes restriction or projection 
onto $S(t)$ and $Ric_{S(t)}$ is the intrinsic Ricci curvature of $S(t)$. 
(The equation (4.15) may be derived from (3.13) by setting $u = 1$ and 
$\sigma = 0$). We recall, as above, that the bar has been removed from 
the notation. 

  As above, the metrics $g_{0}$ and $g_{1}$ are assumed to have a fixed 
geodesic defining function $t$ with common boundary metric $\gamma$ and 
commong geodesic boundary coordinates. Taking the difference of the 
equation (4.15) evaluated on $g_{1}$ and $g_{0}$ gives the following 
equation for $k = g_{1} - g_{0}$ as in (4.8):
\begin{equation} \label{e4.16}
t\ddot k - (n-1)\dot k = tr(\dot k) g_{0}^{T} + 2t(Ric_{S(t)}^{1} - Ric_{S(t)}^{0}) 
+ O(t)k + O(t^{2})\dot k,
\end{equation}
where $O(t^{k})$ denotes terms of order $t^{k}$ with coefficients depending 
smoothly on $g_{0}$. One has $Ric_{S(t)} = D_{x}^{2}(g_{ij})$ is a $2^{\rm nd}$ 
order operator on $g_{ij}$, so that (4.16) gives
\begin{equation} \label{e4.17}
t\ddot k - (n-1)\dot k = tr(\dot k) g_{0}^{T} + 2tD_{x}^{2}(k) + O(t)k + 
O(t^{2})\dot k,
\end{equation}

   The (positive) indicial root of the trace-free part of (4.16) or (4.17) 
is $n$, in that the formal power series solution of (4.17) has undetermined 
coefficient at order $t^{n}$, as in the Fefferman-Graham expansion (4.3)-(4.4). 
The hypothesis (4.7) implies that 
\begin{equation} \label{e4.18}
k = o(t^{n}), 
\end{equation}
so that this $n^{\rm th}$ order coefficient vanishes. However, taking the 
trace of (4.17) gives
$$t \, tr\ddot k - (2n-1)tr \dot k = tr (O(t)k + O(t^{2}\dot k)) + 
2t \, tr (D_{x}^{2}(k)),$$
which has indicial root $2n$. To see that $tr k$ is in fact formally determined 
at order $2n$, one uses the trace of the Riccati equation (3.10), (with $u = 1$ 
and $\sigma = 0$),  which gives
\begin{equation}\label{e4.19}
\dot H  + |A|^{2} = -Ric(T,T).
\end{equation}
Via (4.9), this is easily seen to be equivalent to 
$$t\dot H - H = -t|A|^{2}.$$
This holds for each compactified metric $g_{1}$ and $g_{0}$, and so taking 
the difference, and computing as in (4.16)-(4.17) gives the equation
\begin{equation} \label{e4.20}
t\frac{d^{2}}{dt^{2}}(tr k)  - \frac{d}{dt}(tr k) = O(t)k + O(t^{2})\dot k .
\end{equation}

  The positive indicial root of (4.20) is 2, and by (4.7), the $O(t^{2})$ 
component of the formal expansion of $tr k$ vanishes. Similarly, the trace-free 
part $k_{0}$ of $k$ satisfies the equation
\begin{equation} \label{e4.21}
t\ddot k_{0} -(n-1)\dot k_{0} = 2t(D_{x}^{2}(k))_{0} + [O(t)k]_{0} + 
[O(t^{2})\dot k]_{0} , 
\end{equation}
with indicial root $n$. As in [18], by repeated differentiation of (4.20) 
and (4.21) it follows from (4.7) that the formal expansion of $k$ vanishes. 

  Next we show that (4.8) holds weakly. 

\begin{lemma} \label{l4.2}
Suppose $k = o(t^{n})$ weakly, in that, with respect to the compactified metric 
$(S(t), g)$, ($g = g_{0}$), 
\begin{equation} \label{e4.22}
\int_{S(t)}\langle k, \phi \rangle = o(t^{n}), \ {\rm as} \ t \rightarrow 0 ,
\end{equation}
where $\phi$ is any symmetric bilinear form, $C^{\infty}$ smooth up to 
$U = \partial \Omega$ and vanishing to infinite order on $\partial'\Omega 
= \bar \Omega \setminus U$. Then 
\begin{equation} \label{e4.23}
k = o(t^{\nu}), \ {\rm weakly},
\end{equation}
for any $\nu < \infty$, i.e.~{\rm (4.22)} holds, with $\nu$ in place of $n$. 
\end{lemma}

{\bf Proof:} 
Here smoothness is measured with respect to the given geodesic coordinates 
$(t, x^{i})$ covering $\Omega$. The proof proceeds by induction, starting at 
the initial level $n$. As above, the trace-free and pure trace cases are treated 
separately, and so we assume in the following first that $\phi$ is trace-free. 
Pair $k$ with $\phi$ and integrate (4.17) over the level sets $S(t)$ to obtain
\begin{equation} \label{e4.24}
t\int_{S(t)}\langle \ddot k, \phi \rangle - (n-1)\int_{S(t)}\langle \dot k, 
\phi \rangle  = t\int_{S(t)}\langle k, P_{2}(\phi) \rangle + 
\int_{S(t)}\langle O(t)k, \phi \rangle + \int_{S(t)}\langle O(t^{2})\dot k, \phi \rangle .
\end{equation}
Here $P_{2}(\phi)$ is obtained by integrating the $D_{x}^{2}$ term on the right in 
(4.17) by parts over $S(t)$. Thus $P_{2}(\phi)$, and more generally, $P_{k}(\phi)$ 
denote differential operators of order $k$ on $\phi$ with coefficients depending 
on $g$ and $g_{1}$ and their derivatives up to order 2 and so at least continuous 
up to $\bar \partial \Omega$. We use these expressions generically, so their 
exact form may change from line-to-line below. Note also there are no boundary 
terms at $\partial S(t)$ arising from the integration by parts, by the vanishing 
hypothesis on $\partial'\Omega$. 

    For the terms on the right in (4.24) one then has
$$\int_{S(t)}\langle O(t)k, \phi \rangle = t\int_{S(t)}\langle k, 
P_{0}(\phi) \rangle ,$$
while, since $A = O(t)$ and $H = O(t)$, 
$$\int_{S(t)}\langle O(t^{2})\dot k, \phi \rangle = 
t^{2}\int_{S(t)}\langle \dot k, P_{0}(\phi) \rangle 
= t^{2}\frac{d}{dt}\int_{S(t)}\langle k, P_{0}(\phi) \rangle - 
t^{2}\int_{S(t)}\langle k, P_{1}(\phi) \rangle .$$
Similarly, for the terms on the left in (4.24), one has
$$\int_{S(t)}\langle \dot k, \phi \rangle  = 
\frac{d}{dt}\int_{S(t)}\langle k, \phi \rangle  - t \int_{S(t)}\langle k, P_{1}(\phi) \rangle ,$$
while
$$\int_{S(t)}\langle \ddot k, \phi \rangle = 
\frac{d^{2}}{dt^{2}}\int_{S(t)}\langle k, \phi \rangle 
- 2t\frac{d}{dt}\int_{S(t)}\langle k, P_{1}(\phi) \rangle + 
\int_{S(t)}\langle k, P_{1}(\phi) \rangle  + t\int_{S(t)}\langle k, P_{2}(\phi) \rangle .$$

   Now let 
$$f = f(t) = \int_{S(t)}\langle k, \phi \rangle  .$$
Then the computations above give
\begin{equation} \label{e4.25}
t\ddot f - (n-1)\dot f = t\int_{S(t)}\langle k, P_{2}(\phi) \rangle + (1 + t^{2})\int_{S(t)}
\langle  k, P_{1}(\phi) \rangle
\end{equation}
$$+ \frac{d}{dt}\int_{S(t)}t^{2}\langle k, P_{0}(\phi) \rangle 
+ \frac{d}{dt}\int_{S(t)}t\langle k, P_{1}(\phi) \rangle .$$

   First observe that 
\begin{equation}\label{e4.26}
\int_{S(t)}\langle k, \phi \rangle  = o(t^{n}) \Rightarrow \int_{S(t)}\langle k, 
P_{k}(\phi) \rangle = o(t^{n}),
\end{equation}
for all $C^{\infty}$ forms $\phi$ vanishing to infinite order at $\partial'\Omega$. 
For if the left side of (4.26) holds, then $\int_{S(t)}\langle k, \partial^{k}\phi 
\rangle  = o(t^{n})$, since the hypotheses on $\phi$ are closed under differentiation. 
The coefficients of $P_{k}$ are at least continuous, and it is elementary to verify 
that if $\int_{S(t)}\langle k, \partial^{k}\phi \rangle = o(t^{n})$, then 
$\int_{S(t)}\langle k, \phi \partial^{k}\phi \rangle = o(t^{n})$, for any function 
$\phi$ continuous on $\bar \Omega$. Note that the same result holds with $p$ in 
place of $n$, for any $p < \infty$. 

  It follows from (4.26) and the initial hypothesis (4.22) that the first two terms 
on the right in (4.25) are $o(t^{n})$ as $t \rightarrow 0$. Since $t\ddot f - 
(n-1)\dot f = t^{n}\frac{d}{dt}(\frac{\dot f}{t^{n-1}})$, this gives
$$\frac{d}{dt}(\frac{\dot f}{t^{n-1}}) = o(1) + t^{-n}\frac{d}{dt}\int_{S(t)}t\langle k, 
P_{1}(\phi) \rangle + t^{-n}\frac{d}{dt}\int_{S(t)}t^{2}\langle k, P_{0}(\phi) \rangle .$$
Integrating from $0$ to $t$ implies
$$\frac{\dot f}{t^{n-1}} = o(t) + t^{-n+1}\int_{S(t)}\langle k, P_{1}(\phi) \rangle 
+ n\int_{0}^{t}t^{-n}\int_{S(t)}\langle k, P_{1}(\phi) \rangle + c_{1} = o(t) 
+ c_{1} ,$$
where $c_{1}$ is a constant. A further integration using (4.26) again gives
\begin{equation} \label{e4.27}
f = o(t^{n+1}) + c_{1}'t^{n} + c_{2},
\end{equation}
where $c_{1}' = \frac{c_{1}}{n}$. Once more by (4.22), this implies that
$$f = o(t^{n+1}).$$
Note the special role played by the indicial root $n$ here; if instead one had only 
$k = O(t^{n})$, then the argument above does not give $k = O(t^{n+1})$ weakly. 

  This first estimate holds in fact for any given trace-free $\phi$ which is $C^{2}$ 
on $\bar \Omega$, and vanishing to first order on $\partial'\Omega$. Working in the 
same way with the trace equation (4.20) shows that the same result holds for pure 
trace terms. In particular, it follows that
\begin{equation} \label{e4.28}
k = o(t^{n+1}) \ {\rm weakly} .
\end{equation}

  One now just repeats this argument inductively, with the improved estimate 
(4.28) in place of (4.22), using (4.26) inductively. Note that each inductive 
step requires higher differentiability of the test function $\phi$ and its higher 
order vanishing at $\partial'\Omega$. 

{\endproof}

  Lemma 4.2 proves that $k = k_{\alpha\beta} = O(t^{\nu})$ weakly, for any 
$\nu < \infty$. As discussed in \S 3, the transition from geodesic boundary 
coordinates to $H$-harmonic coordinates is $C^{2,\alpha}$ and hence
\begin{equation}\label{e4.29}
h = h_{\alpha\beta} = O(t^{\nu}),
\end{equation}
weakly, with the level sets $S(t)$ replaced by $\Sigma_{\tau}$. Next, as in 
Remark 3.5 and the proof of Theorem 3.1, the equations (4.11)-(4.13) satisfy 
elliptic estimates, and elliptic regularity in weighted H\"older spaces, cf.~[26], [20], 
shows that the weak decay (4.29) implies strong or pointwise decay, i.e.~(4.14) holds. 
The proof of Theorem 4.1 and thus Theorem 1.2 is now completed as before in the 
$C^{\infty}$ smooth case. 

{\endproof}

\begin{remark} \label{r 4.3}
{\rm In [3, Thm.~3.2], a proof of unique continuation of conformally compact 
Einstein metrics was given in dimension 4, using the fact that the 
compactified metric $\widetilde g$ in (1.1) satisfies the Bach equation, 
together with the Calder\'on uniqueness theorem. However, the proof 
in [3] used harmonic coordinates; as discussed in \S 2, such coordinates 
do not preserve the Cauchy data. I am grateful to Robin Graham for pointing 
this out. Theorem 1.2 thus corrects this error, and generalizes the result to any 
dimension. }
\end{remark}

  For the work to follow in \S 5, we note that Theorem 4.1 also holds for 
linearizations of the Einstein equations, i.e.~forms $k$ satisfying 
\begin{equation}\label{e4.30}
\frac{d}{dt}(Ric_{g+tk} + n(g+tk))|_{t=0} = 0. 
\end{equation}
Thus, if $k$ satisfies (4.30) and the analog of (4.7), i.e.~$|k| = o(t^{n})$, then 
$k$ is pure gauge in $\Omega$, in that $k = \delta^{*}Z$, where $Z$ is a vector 
field on $\Omega$ with $Z = 0$ on $\partial \Omega$. The proof of this is exactly 
the same as the proof of Theorem 4.1, replacing the finite difference $k = g_{1} - 
g_{0}$ by an infinitesimal difference. 

  This has the following consequence:
\begin{corollary}\label{c4.4}
Let $(M, g)$ be a conformally compact Einstein manifold with metric $g$ having a 
$C^{2,\alpha}$ geodesic compactification. Suppose the topological condition 
{\rm (1.5)} holds, i.e.~$\pi_{1}(M, \partial M) = 0$. 

  If $k$ is an infinitesimal Einstein deformation on $M$ as in {\rm (4.30)}, 
in divergence-free gauge, i.e.
\begin{equation}\label{e4.31}
\delta k = 0,
\end{equation}
with $k = o(t^{n})$ on approach to $\partial M$, then
$$k = 0 \ \ {\rm on} \ \ M.$$
\end{corollary}

{\bf Proof:}  The topological condition (1.5), together with the same analytic 
continuation argument at the end of the proof of Theorem 3.1, (cf.~also (3.35)), 
implies that $k$ is pure gauge globally on $M$, in that $k = \delta^{*}Z$ on 
$M$ with $Z = 0$ on $\partial M$. (Recall that (1.5) implies that $\partial M$ 
is connected). From (4.31), one then has
$$\delta \delta^{*}Z = 0,$$
on $M$. Pairing this with $Z$ and integrating over $B(t)$, it follows that 
$$\int_{B(t)}|\delta^{*}Z|^{2} = \int_{S(t)}\delta^{*}Z(Z, N),$$
where $N$ is the unit outward normal. Since $|Z|_{g}$ is bounded and 
$|\delta^{*}Z|vol(S(t)) = o(1)$, (since $|k| = o(t^{n})$), it follows that 
$$\int_{M}|\delta^{*}Z|^{2} = 0,$$
which gives the result.

{\endproof}

  Of course, analogs of these results also hold for bounded domains, via the 
proof of Theorem 3.1; the verification is left to the reader. 

\begin{remark} \label{r4.5}
{\rm The analogue of Proposition 3.6 most likely also holds in the setting of 
conformally compact metrics, for fields $\tau$ whose Euler-Lagrange 
equation is a diagonal system of Laplace-type operators to leading 
order, as in (3.37) or (3.38). The proof of this is basically the same 
as that of Proposition 3.6, using the proof of Theorem 1.2 and with 
the Mazzeo unique continuation result in place of that of Calder\'on. 
However, we will not carry out the details of the proof here. }
\end{remark}

\section{Isometry Extension and the Constraint Equations.}
\setcounter{equation}{0}

  In this section, we prove Theorem 1.3 that continuous groups of isometries at the 
boundary extend to isometries in the interior of complete conformally compact Einstein 
metrics and relate this issue in general to the constraint equations induced by the 
Gauss-Codazzi equations.  

  We begin with the following elementary consequence of Theorem 4.1. 

\begin{proposition} \label{p5.1}
Let $(\Omega, g)$ be a $C^{n}$ polyhomogeneous conformally compact Einstein metric 
on a domain $\Omega \simeq B^{n+1}$ with boundary metric $\gamma$ on $\partial 
\Omega \simeq B^{n}$. Suppose $X$ is a Killing field on $(\partial \Omega, \gamma)$ 
and 
\begin{equation} \label{e5.1}
{\mathcal L}_{X}g_{(n)} = 0,
\end{equation}
where $g_{(n)}$ is the $n^{\rm th}$ term in the Fefferman-Graham expansion {\rm (4.3)} 
or {\rm (4.4)}.

  Then $X$ extends to a Killing field on $(\Omega, g)$. 
\end{proposition}

{\bf Proof:} 
Extend $X$ to a smooth vector field on $\Omega$ by requiring $[X, N] = 0$, where 
$N = \nabla \log t$ and $t$ is the geodesic defining function determined by $g$ and 
$\gamma$. Let $\phi_{s}$ be the corresponding 1-parameter group of diffeomorphisms 
and set $g_{s} = \phi^{*}_{s}g$. Then $t$ is the geodesic defining function for 
$g_{s}$ for any $s$, and the pair $(g, g_{s})$ satisfy the hypotheses of Theorem 4.1. 
Theorem 4.1 then implies that $g_{s}$ is isometric to $g$, i.e.~there 
exist diffeomorphisms $\psi_{s}$ of $\Omega$, equal to the identity on 
$\partial \Omega$, such that $\psi_{s}^{*}\phi_{s}^{*}g = g$. Thus 
$\phi_{s}\circ \psi_{s}$ is a 1-parameter group of isometries of $g$ defined in 
$\Omega$, with $Y$ the corresponding Killing field. (In fact, $Y = X$, since any Killing 
field $Y$ tangent to $\partial\Omega$ preserves the geodesics tangent to $N$, and so 
$[Y, N] = 0$. This determines $Y$ uniquely in terms of its value at $\partial \Omega$. 
Since $X$ satsifies the same equation with the same initial value, this gives the claim). 

{\endproof}

  We point out that the the same result, and proof, also hold in the case of 
Einstein metrics on bounded domains, via Theorem 3.1. The condition (5.1) is 
of course replaced by ${\mathcal L}_{X}A = 0$. For some examples and discussion 
in the bounded domain case, see [1], [2]. 

  Suppose now that $(M, g)$ is a (global) conformally compact Einstein metric and 
there is a domain $\Omega$ as in Proposition 5.1 contained in $M$ on which (5.1) 
holds. Then by analytic continuation as discussed at the end of the proof of Theorem 
3.1, $X$ extends to a local Killing field on all of $M$, i.e.~$X$ extends to a Killing 
field on the universal cover $\widetilde M$. In particular, if the condition (3.35) holds, 
i.e.
$$\pi_{1}(\partial M) \rightarrow \pi_{1}(M) \rightarrow 0,$$
then $X$ extends to a global Killing field on $M$. Again, the same result holds 
in the context of bounded domains. 

\begin{remark}\label{r5.2}
{\rm A natural analogue of Proposition 5.1 holds for conformal Killing fields on 
$(\partial \Omega, \gamma)$, i.e.~vector fields which preserve the conformal class 
$[\gamma]$ at conformal infinity. Such vector fields satisfy the conformal Killing 
equation
\begin{equation}\label{e5.2}
\hat{\mathcal L}_{X}\gamma = {\mathcal L}_{X}\gamma - \frac{tr(\mathcal L_{X}\gamma)}{n}
\gamma = 0 .
\end{equation}
Namely, since we are working locally, it is well-known - and easy to prove - that any 
non-vanishing conformal Killing field is Killing with respect to a conformally related 
metric $\widetilde \gamma = \lambda^{2}\gamma$, so that 
$${\mathcal L}_{X}\widetilde \gamma = 0 .$$
Hence, if ${\mathcal L}_{X}\widetilde g_{(n)} = 0$, then Proposition 5.1 implies that 
$X$ extends to a Killing field on $\Omega$. 

   One may express $\widetilde g_{(n)}$ in terms of $\lambda$ and the lower order 
terms $g_{(k)}$, $k < n$ in the Fefferman-Graham expansion (4.3)-(4.4); however, the 
expressions become very complicated for $n$ even and large, cf.~[17]. Thus, while the 
equation (5.2) is conformally invariant, the corresponding conformally invariant 
equation for $g_{(n)}$ will be complicated in general. }
\end{remark}

   Next we consider the constraint equations (4.5) in detail, i.e.
\begin{equation}\label{e5.3}
\delta \tau_{(n)} = 0 \ \ {\rm and} \ \ tr \,\tau_{(n)} = a_{(n)},
\end{equation}
where $\tau_{(n)} = g_{(n)} + r_{(n)}$; $r_{(n)}$ and $a_{(n)}$ are explicitly determined 
by the boundary metric $\gamma = g_{(0)}$ and its derivatives up to order $n$. Both vanish 
when $n$ is odd. 

  As will be seen below, the most important issue is the divergence constraint in 
(5.3), which arises from the Gauss-Codazzi equations. To see this, in the setting 
of \S 4, on $S(t) \subset (M, g)$, the Gauss-Codazzi equations are
\begin{equation} \label{e5.4}
\delta(A - Hg) = -Ric(N, \cdot),
\end{equation}
as 1-forms on $S(t)$; here $N = -t\partial_{t}$ is the unit outward normal. The same 
equation holds on a geodesic compactification $(M, \bar g)$. If $g$ is Einstein, then 
$Ric(N, \cdot) = \bar Ric(\bar N, \cdot) = 0$; the latter equality follows from (4.9). 
The equation (5.4) holds for all $t$ small, and differentiating $(n-1)$ times with 
respect to $t$ gives rise to the divergence constraint in (5.3). 

   The Gauss-Codazzi equations are not used in the derivation and properties of 
the Fefferman-Graham expansion (4.3)-(4.4) per se. The derivation of these equations 
involves only the tangential $(ij)$ part of the Ricci curvature. The asymptotic 
behavior of the normal $(00)$ part of the Ricci curvature gives rise to the trace 
constraint in (5.3), cf.~(4.19)-(4.20). 

\medskip

  Let ${\mathcal T}$ be the space of pairs $(g_{(0)}, \tau_{(n)})$ satisfying (5.3). 
If $\tau_{(n)}^{0}$ is any fixed solution of (5.3), then any other solution with 
the same $g_{(0)}$ is of the form $\tau_{(n)} = \tau_{(n)}^{0} + \tau$, where 
$\tau$ is transverse-traceless (TT) with respect to $g_{(0)}$. (Of course if $n$ 
is odd, one may take $\tau_{(n)}^{0} = 0$). The space ${\mathcal T}$ is naturally a 
bundle over $Met(\partial M)$ with fiber at $\gamma$ an affine space of symmetric tensors 
and is a subset of the product 
$Met(\partial M) \times {\mathbb S}^{2}(\partial M) \simeq T(Met(\partial M))$. Let 
\begin{equation}\label{e5.5}
\pi: {\mathcal T} \rightarrow Met(\partial M)
\end{equation}
be the projection onto the base space $Met(\partial M)$, (the first factor projection). 

  By the discussion in \S 4, $(g_{(0)}, \tau_{(n)}) \in  {\mathcal T}$ if and only 
if the corresponding pair $(g_{(0)}, g_{(n)})$ determine a formal polyhomogenous solution 
to the Einstein equations near conformal infinity, i.e.~formal series solutions 
containing $\log$ terms, as in (4.3)-(4.4). In fact, if  $g_{(0)}$ and $g_{(n)}$ are 
real-analytic on $\partial M$, a result of Kichenassamy [24] implies that the series 
(4.3) or (4.4) converges, and gives an Einstein metric $g$, defined in a neighborhood 
of $\partial M$. The metric $g$ is complete near $\partial M$ and has a conformal 
compactification inducing the given data $(g_{(0)}, g_{(n)})$ on $\partial M$. 
Here we recall from the discussion in \S 4 that all coefficients of the expansion 
(4.3) or (4.4) are determined by $g_{(0)}$ and $g_{(n)}$. 

  In this regard, consider the following:

  {\bf  Problem.} Is $\pi: {\mathcal T} \rightarrow Met(\partial M)$ an open map? 
Thus, given any $(g_{(0)}, \tau_{(n)}) \in {\mathcal T}$ and any boundary metric 
$\widetilde g_{(0)}$ sufficiently close to $g_{(0)}$, does there exist 
$\widetilde \tau_{(n)}$ close to $\tau_{(n)}$ such that $(\widetilde g_{(0)}, 
\widetilde \tau_{(n)}) \in {\mathcal T}$.

  Although $\pi$ is obviously globally surjective, the problem above is whether 
$\pi$ is locally surjective. For example, a simple fold map $x \rightarrow x^{3}-x$ 
is not locally surjective near $\pm\sqrt{3}/3$. Observe that the trace condition in 
(5.3) imposes no constraint on $g_{(0)}$; given any $g_{(0)}$, it is easy to find 
$g_{(n)}$ such that $tr_{g_{(0)}}(g_{(n)}+r_{(n)}) = a_{(n)}$; this equation can readily 
be solved algebraically for many $g_{(n)}$. 

  By the inverse function theorem, it suffices, (and is probably also necessary), to 
examine the problem above at the linearized level. However the linearization of the 
divergence condition in (5.3) gives a non-trivial constraint on the variation $h_{(0)}$ 
of $g_{(0)}$. Namely, the linearization in this case gives 
\begin{equation} \label{e5.6} 
\delta'(\tau_{(n)}) + \delta(\tau_{(n)})' = 0 ,
\end{equation}
where $\delta ' = \frac{d}{du}\delta_{g_{(0)}+uh_{(0)}}$, and similarly for 
$(\tau_{(n)})'$.

  Whether (5.6) is solvable for any $h_{(0)} \in S^{2}(\partial M)$ 
depends on the data $g_{(0)}$ and $g_{(n)}$. For example, it is trivially solvable when 
$\tau_{(n)} = 0$. For compact $\partial M$, one has 
\begin{equation} \label{e5.7}
\Omega^{1}(\partial M) = Im \delta \oplus Ker \delta^{*}, 
\end{equation}
where $\Omega^{1}$ is the space of 1-forms, so that solvability in general requires that 
\begin{equation}\label{e5.8}
\delta'(\tau_{(n)}) \in Im \delta = (Ker \delta^{*})^{\perp}.
\end{equation}
Of course $Ker \delta^{*}$ is exactly the space of Killing fields on $(\partial M, \gamma)$, 
and so this space serves as a potential obstruction space.  

  Clearly then $\pi$ is locally surjective when $(\partial M, g_{(0)})$ has no Killing fields. 
On the other hand, it is easy to construct examples where $(\partial M, \gamma)$ does 
have Killing fields and $\pi$ is not locally surjective: 

\begin{example} \label{ex5.3}
{\rm  Let $(\partial M, g_{(0)})$ be the flat metric on the $n$-torus $T^{n}$, $n \geq 3$, 
and define $g_{(n)} = -(n-2)(d\theta^{2})^{2} +(d\theta^{3})^{2} + \cdots + (d\theta^{n})^{2}$. 
Then $g_{(n)}$ is transverse-traceless with respect to $g_{(0)}$. Let $f = f(\theta^{1})$. 
Then $\hat g_{(n)} = fg_{(n)}$ is still TT with respect to $g_{(0)}$, so that 
$(g_{(0)}, \hat g_{(n)}) \in {\mathcal T}$, at least for $n$ odd. 

   It is then not difficult to see via a direct calculation, or more easily via 
Proposition 5.4 below, that (5.8) does not hold, so that $\pi$ is not locally surjective. }
\end{example}

   Next we relate these two issues, i.e.~the general solvability of the divergence 
constraint (5.8) and the extension of Killing fields on the boundary into the bulk. 
The following result holds for general $\phi \in S^{2}(\partial M)$ with $\delta \phi = 0$, 
but will only be used in the case $\phi = \tau_{(n)}$ on $\partial M$. 

\begin{proposition} \label{p5.4}
If $X$ is a Killing field on $(\partial M, \gamma)$, with $\partial M$ compact, then 
\begin{equation}\label{e5.9}
\int_{\partial M}\langle {\mathcal L}_{X}\tau_{(n)}, h_{(0)} \rangle dV= 
-2\int_{\partial M}\langle \delta'(\tau_{(n)}), X \rangle dV,
\end{equation}
where $\delta' = \frac{d}{ds}\delta_{\gamma + sh_{(0)}}$. 
In particular, {\rm (5.1)} holds for all Killing fields on $(\partial M, \gamma)$ if 
and only if the linearized divergence constraint vanishes, i.e.~{\rm (5.8)} holds. 
\end{proposition}

{\bf Proof:}
We will carry out the computation in general, and only use the condition that $X$ 
is a Killing field at the end of the proof. To simplify the notation, set 
$h_{(0)} = h$.  

  The following formulas will be used here and in work to follow:
\begin{equation} \label{e5.10}
{\mathcal L}_{V}\phi = \nabla_{V}\phi + 2\nabla V\circ \phi,
\end{equation}
which is standard, and 
\begin{equation}\label{e5.11}
(\delta^{*})'X = \frac{1}{2}\nabla_{X}h + \delta^{*}X\circ h,
\end{equation}
cf.~[11] for example. Here $\phi \circ \psi$ is the symmetrized product; in an 
orthonormal frame, $(\phi \circ \psi)(e_{i}, e_{j}) = \frac{1}{2}(\langle \phi(e_{i}), 
\psi(e_{j}) \rangle + \langle \phi(e_{j}), \psi(e_{i}) \rangle)$. 

  To begin $\int_{\partial M}\langle {\mathcal L}_{X}\tau_{(n)}, h \rangle = 
\int_{\partial M}\langle \nabla_{X}\tau_{(n)}, h \rangle + 
2\langle \nabla X \circ \tau_{(n)}, h \rangle$. Since $h$ is symmetric, 
$\langle \nabla X \circ \tau_{(n)}, h \rangle = \langle \delta^{*}X 
\circ \tau_{(n)}, h \rangle$. For the first term, write 
$\langle \nabla_{X}\tau_{(n)}, h \rangle = X\langle \tau_{(n)}, h \rangle 
- \langle \tau_{(n)}, \nabla_{X}h \rangle$. The first term here integrates to 
$\delta X \langle \tau_{(n)}, h \rangle$, while by (5.11), the second term is 
$- \langle \tau_{(n)}, \nabla_{X}h \rangle = -2 \langle \tau_{(n)}, 
(\delta^{*})'X \rangle + 2\langle \tau_{(n)}, \delta^{*}X\circ h \rangle$. 

  Next, a straightforward computation using the fact that $\delta \tau_{(n)} = 0$ 
gives 
$$\int_{\partial M}\langle \tau_{(n)}, (\delta^{*})'X \rangle dV = 
\int_{\partial M}\langle (\delta')(\tau_{(n)}), X \rangle dV$$
$$+ 2\int_{\partial M}\langle \tau_{(n)}\circ \delta^{*}X,  h \rangle - 
\frac{1}{2}\int_{\partial M}\langle \tau_{(n)},\delta^{*}X \rangle tr h dV .$$ 
The last two terms come from variation of the metric and volume form. 
Combining these computations gives 

\begin{equation}\label{e5.12}
\int_{\partial M}\langle {\mathcal L}_{X}\tau_{(n)}, h \rangle = 
-2\int_{\partial M}\langle \delta'(\tau_{(n)}), X \rangle 
+ \int_{\partial M}[\delta X \langle \tau_{(n)}, h \rangle + {\tfrac{1}{2}}
\langle \tau_{(n)}, \delta^{*}X \rangle tr h]dV .
\end{equation}

This gives (5.9) when $X$ is Killing, i.e.~$\delta^{*}X = 0$; note that 
${\mathcal L}_{X}r_{(n)} = 0$ in this case, since $r_{(n)}$ is determined by 
the boundary metric. 

  To prove the last statement, by (5.9), (5.1) holds if and only if 
$\int_{\partial M}\langle \delta'(\tau_{(n)}), X \rangle = 0$, for all 
variations $h$. If (5.8) holds, then $\delta'(\tau_{(n)}) = \delta h_{(n)}'$, 
for some $h_{(n)}'$ and so $\int_{\partial M}\langle \delta'(\tau_{(n)}), X \rangle 
= \int_{\partial M}\langle h_{(n)}', \delta^{*}X \rangle = 0$, since $X$ is 
Killing. The converse of this argument holds equally well. 

{\endproof}

  Proposition 5.4 implies that in general, Killing fields on $\partial M$ do not 
extend to Killing fields in a neighborhood of $\partial M$, (cf.~Example 5.3). 
(Exactly the same result and proof hold in the bounded domain case, when the term 
$\tau_{(n)}$ is replaced by $A - Hg$). 

  Now as noted above, whether isometry extension holds or not depends on the 
term $\tau_{(n)} = g_{(n)}+r_{(n)}$, or more precisely on the relation of the 
boundary metric $g_{(0)}$ with $\tau_{(n)}$. For Einstein metrics which are 
globally conformally compact, the term $\tau_{(n)}$ is determined, up to a finite 
dimensional moduli space, by the boundary metric $g_{(0)}$; (this is discussed 
further below). Thus, whether isometry extension holds or not is quite a delicate 
issue; if so, it must depend crucially on the global structure of $(M, g)$. 

\medskip

  Before beginning the proof of Theorem 1.3, we first need to discuss some 
background material from [5]-[6].

  Let $E_{AH}$ be the space of conformally compact, or equivalently asymptotically 
hyperbolic Einstein metrics on $M$ which have a $C^{\infty}$ polyhomogeneous
conformal compactification with respect to a fixed smooth defining function $\rho$, 
as in (1.1). In [5], it is shown that $E_{AH}$ is a smooth, infinite dimensional 
manifold. One has a natural smooth boundary map 
\begin{equation}\label{e5.13}
\Pi: E_{AH} \rightarrow Met(\partial M),
\end{equation}
sending $g$ to its boundary metric $\gamma$. 

  The moduli space ${\mathcal E}_{AH}$ is the quotient $E_{AH}/{\mathcal D}_{1}$, 
where ${\mathcal D}_{1}$ is the group of smooth (polyhomogeneous) diffeomorphisms 
$\phi$ of $M$ equal to the identity on $\partial M$. Thus, $g' \sim g$ if 
$g' = \phi^{*}g$, with $\phi \in {\mathcal D}_{1}$. Changing the defining function 
$\rho$ in (1.1) changes the boundary metric conformally. Also, if 
$\phi \in {\mathcal D}_{1}$ then $\rho \circ \phi$ is another defining function, 
and all defining functions are of this form near $\partial M$. Hence if 
${\mathcal C}$ denotes the space of smooth conformal classes of metrics on 
$\partial M$, then the boundary map (5.13) descends to a smooth map 
\begin{equation}\label{e5.14}
\Pi: {\mathcal E}_{AH} \rightarrow {\mathcal C}
\end{equation}
independent of the defining function $\rho$. Either boundary map $\Pi$ in (5.13) 
or (5.14) is smooth and Fredholm, of Fredholm index 0. 

 The linearization of the Einstein operator $Ric_{g} + ng$ at an Einstein metric 
$g$ is given by 
\begin{equation}\label{e5.15}
\hat L = (Ric_{g} + ng)' = \tfrac{1}{2}D^{*}D -  R - \delta^{*}\beta, 
\end{equation}
acting on the space of symmetric 2-tensors $S^{2}(M)$ on $M$, cf.~[11]. Here, (as 
in \S 3), $\beta$ is the Bianchi operator, $\beta(h) = \delta h + \frac{1}{2}d tr h$, 
Thus, $h \in T_{g}E_{AH}$ if and only if 
$$\hat L(h) = 0.$$
The operator $\hat L$ is not elliptic, due to the $\delta^{*}\beta$ term. As is 
well-known, this arises from the diffeomorphism group, and to obtain an elliptic 
linearization, one needs a gauge choice to break the diffeomorphism invariance 
of the Einstein equations. We will use a slight modification of the Bianchi 
gauge introduced in [12]. 

  To describe this, given any fixed $g_{0} \in E_{AH}$ with geodesic defining 
function $t$ and boundary metric $\gamma_{0}$, let $\gamma$ be a boundary metric 
near $\gamma_{0}$ and define the hyperbolic cone metric $g_{\gamma}$ on $\gamma$ 
by setting 
$$g_{\gamma} = t^{-2}(dt^{2} + \gamma);$$
$g_{\gamma}$ is defined in a neighborhood of $\partial M$. Next, set
\begin{equation}\label{e5.16}
g(\gamma) = g_{0} + \eta(g_{\gamma} - g_{\gamma_{0}}),
\end{equation}
where $\eta$ is a non-negative cutoff function supported near $\partial M$ with 
$\eta = 1$ in a small neighborhood of $\partial M$. Any conformally compact metric 
$g$ near $g_{0}$, with boundary metric $\gamma$ then has the form 
\begin{equation}\label{e5.17}
g = g(\gamma) + h,
\end{equation}
where $|h|_{g_{0}} = O(t^{2})$; equivalently $\bar h = t^{2}h$ satisfies $\bar h_{ij} 
= O(t^{2})$ in any smooth coordinate chart near $\partial M$. The space of such 
symmetric bilinear forms $h$ is denoted by ${\mathbb S}_{2}(M)$ and the space of 
metrics $g$ of the form (5.17) is denoted by $Met_{AH}$. 

  The Bianchi-gauged Einstein operator, (with background metric $g_{0}$), is defined 
by
\begin{equation}\label{e5.18}
\Phi_{g_{0}}: Met_{AH} \rightarrow {\mathbb S}_{2}(M)
\end{equation}
$$\Phi_{g_{0}}(g) = \Phi (g(\gamma) + h) = Ric_{g} + ng + 
(\delta_{g})^{*}\beta_{g(\gamma)}(g),$$
where $\beta_{g(\gamma)}$ is the Bianchi operator with respect to $g(\gamma)$. By 
[12, Lemma I.1.4],
\begin{equation}\label{e5.19}
Z_{AH} \equiv  \Phi^{-1}(0)\cap\{Ric <  0\} \subset E_{AH}, 
\end{equation}
where $\{Ric < 0\}$ is the open set of metrics with negative Ricci curvature. 
In fact, if $g \in E_{AH}$ is close to $g_{0}$, and $\Phi (g) = 0$, then 
$\beta_{g(\gamma)}(g) = 0$ and moreover
\begin{equation}\label{e5.20}
\delta_{g(\gamma)}(g) = 0, \ {\rm and} \ tr_{g(\gamma)}(g) = 0. 
\end{equation}
The space $Z_{AH}$ is a local slice for the action of ${\mathcal D}_{1}$ on 
$E_{AH}$: for any $g\in E_{AH}$ near $g_{0}$, there exists a diffeomorphism 
$\phi \in {\mathcal D}_{1}$ such that $\phi^{*}g\in Z_{AH}$, cf.~again [12].

   The linearization of $\Phi$ at $g_{0} \in E_{AH}$ with respect to the $2^{\rm nd}$ 
variable $h$ has the simple form
\begin{equation}\label{e5.21}
(D_{2}\Phi)_{g_{0}}(\dot h) = \tfrac{1}{2}D^{*}D \dot h -  R_{g_{0}}(\dot h), 
\end{equation}
while the variation of $\Phi$ at $g_{0}$ with respect to the $1^{\rm st}$ variable 
$g(\gamma)$ has the form 
\begin{equation}\label{e5.22}
(D_{1}\Phi)_{g_{0}}(\dot g(\gamma)) = (D_{2}\Phi)_{g_{0}}(\dot g(\gamma)) 
- \delta_{g_{0}}^{*}\beta_{g_{0}}(\dot g(\gamma)) = (Ric_{g} + ng)'(\dot g(\gamma)),
\end{equation}
as in (5.15). Clearly $\dot g(\gamma) = \eta t^{-2}\dot \gamma$. 
  The kernel of the elliptic self-adjoint linear operator
\begin{equation}\label{e5.23} 
L = \tfrac{1}{2}D^{*}D -  R 
\end{equation}
acting on the $2^{\rm nd}$ variable $h$, represents the space of non-trivial 
infinitesimal Einstein deformations vanishing on $\partial M$. Let $K$ denote the 
$L^{2}$ kernel of $L$. This is the same as the kernel of $L$ on ${\mathbb S}_{2}(M)$, 
cf.~[12], [26]. An Einstein metric $g_{0} \in E_{AH}$ is called {\it non-degenerate} 
if 
\begin{equation} \label{e5.24}
K = 0. 
\end{equation}

  For $g_{0} \in {\mathcal E}_{AH}$ the kernel $K = K_{g_{0}}$ equals the kernel 
of the linear map $D\Pi: T_{g_{0}}{\mathcal E}_{AH} \rightarrow T_{\Pi(g_{0})}
{\mathcal C}$. Hence, $g_{0}$ is non-degenerate if and only if $g_{0}$ is a regular 
point of the boundary map $\Pi$ in which case $\Pi$ is a local diffeomorphism near 
$g_{0}$. From now on, we denote $g_{0}$ by $g$. 

  The components of any $\kappa \in K$ satisfy the following bounds with respect 
to a geodesic defining function, cf.~[20], [26] or [28,~Prop.~5] for example:
\begin{equation}\label{e5.25}
\kappa = O(t^{n}), \ \ \kappa(N,Y) = O(t^{n+1}), \ \ \kappa (N,N) = 
O(t^{n+1+\mu}),
\end{equation}
where $N = -t\partial_{t}$ is the unit outward normal vector to the $t$-level set 
$S(t)$, $Y$ is any $g$-unit vector tangent to $S(t)$ and $\mu > 0$. Here 
$\kappa = O(t^{n})$ means $|\kappa|_{g} = O(t^{n})$. Also by (5.20), any 
$\kappa \in K$ is transverse-traceless, i.e.
\begin{equation}\label{e5.26}
\delta \kappa = tr \kappa = 0.
\end{equation}

\medskip

  Given this background, we are now ready to begin the proof of Theorem 1.3. 

{\bf Proof of Theorem 1.3.}

\medskip
 
  Let $\bar g = t^{2}g$ be a geodesic compactification of $g$ with 
boundary metric $\gamma$. By the boundary regularity result of [16], 
$\bar g$ is $C^{\infty}$ polyhomogeneous on $\bar M$. It suffices to prove 
Theorem 1.3 for arbitrary 1-parameter subgroups of the isometry group of 
$(\partial M, \gamma)$. Thus, let $\phi_{s}$ be a local 1-parameter group 
of isometries of $\gamma$ with $\phi_{0} = id$, so that
$$\phi_{s}^{*}\gamma  = \gamma .$$ 
The diffeomorphisms $\phi_{s}$ of $\partial M$ may be extended to 
diffeomorphisms of $M$, so that the curve  
\begin{equation} \label{e5.27}
g_{s} = \phi_{s}^{*}g
\end{equation}
is a smooth curve in $E_{AH}$. By construction then, $\Pi[g_{s}] = [\gamma]$, 
so that $[h] = [\frac{dg_{s}}{ds}] \in Ker D\Pi$, for $\Pi$ as in (5.14). One 
may then alter the diffeomorphisms $\phi_{s}$ by composition with diffeomorphisms 
in ${\mathcal D}_{1}$ if necessary, so that $h = \frac{dg_{s}}{ds} \in K_{g}$, 
where $K_{g}$ is the kernel in (5.24). Denoting $h = \kappa$, it follows that 
\begin{equation}\label{e5.28}
\kappa = \delta^{*}X,
\end{equation}
where $X = d\phi_{s}/ds$ is smooth up to $\bar M$. 

  Thus it suffices to prove that $\delta^{*}X = 0$, since this will imply that 
$g_{s} = g$, (when $g_{s}$ is modified by the action of ${\mathcal D}_{1}$). If 
$K_{g} = 0$, i.e.~if $g$ is a regular point of the boundary map $\Pi$, then 
this is now obvious, (from the above), and proves the result in this special 
case; (the proof in this case requires only that $(M, g)$ be $C^{2,\alpha}$ 
conformally compact). 

   To prove the result in general, we will prove that any solution of (5.28) 
necessarily vanishes, so that 
\begin{equation} \label{e5.29}
K \cap Im \delta^{*} = 0.
\end{equation}

  We give two different, (although related), proofs of (5.29), one conceptual and 
one more computational. The first, conceptual, proof involves an understanding of the 
cokernel of the map $D\Pi_{g}$ in $Met(\partial M)$, and so one first needs to give an 
explicit description of this cokernel. To begin, recall the derivative
\begin{equation} \label{e5.30}
(D\Phi)_{g}: T_{g}Met_{AH}(M) \rightarrow  T_{\Phi(g)}{\mathbb S}_{2}(M). 
\end{equation}
Via (5.17), one has $T_{g}Met_{AH} = T_{\gamma}Met(\partial M)\oplus T_{h}{\mathbb S}_{2}(M)$ 
and the derivative with respect to the second factor is given by (5.21). If $K = 0$, then 
$D_{2}\Phi$ is surjective at $g$, (since $D_{2}\Phi$ has index 0), and hence so is $D\Pi$. 
In general, to understand $Coker D\Pi$, we show that $D\Phi$ is always surjective; this 
follows from the claim that for any non-zero $\kappa\in K$ there is a tangent vector 
$\dot g(\gamma) \in T_{\gamma}Met(\partial M) \subset T_{g}Met_{AH}$ such that
\begin{equation} \label{e5.31}
\int_{M}\langle (D_{1}\Phi)_{g}(\dot g(\gamma)), \kappa \rangle dV_{g} \neq 0. 
\end{equation}
Thus, the boundary variations $\dot g(\gamma)$ satisfying (5.31) for some $\kappa$ 
correspond to the cokernel. To prove (5.31), let $B(t) = \{x \in M: t(x) \geq t\}$ 
and $S(t) = \partial B(t) = \{x \in M: t(x) = t\}$. Apply the divergence theorem to 
the integral (5.31) over $B(t)$; twice for the Laplace 
term in (5.22) and once for the $\delta^{*}$ term in (5.22). Since 
$$\kappa \in Ker L \ {\rm and} \  \delta \kappa = 0,$$
it follows that the integral (5.31) reduces to an integral over the boundary, 
and gives
\begin{equation} \label{e5.32}
\int_{B(t)}\langle (D_{1}\Phi)_{g}(\dot g(\gamma), \kappa \rangle dV_{g} = 
{\tfrac{1}{2}}\int_{S(t)}(\langle \dot g(\gamma), \nabla_{N}\kappa \rangle  -  
\langle \nabla_{N}\dot g(\gamma), \kappa \rangle  - 
2\langle \beta(\dot g(\gamma)), \kappa(N) \rangle )dV_{S(t)}. 
\end{equation}
Of course $dV_{S(t)} = t^{-n}dV_{\gamma} + O(t^{-(n-1)})$. By (5.25) the last 
term in (5.32) is then $O(t)$ and so may be ignored. Let
\begin{equation}\label{e5.33}
\widetilde \kappa = t^{-n}\kappa,
\end{equation}
so that by (5.25), $|\widetilde \kappa|_{g} \leq C$. From the definition (5.16), a 
straightforward computation shows that near $\partial M$, 
$$\dot g(\gamma) = t^{-2}\dot \gamma, \ \ {\rm and} \ \ \nabla_{N}\dot g(\gamma) = 0.$$ 
Note that $|\dot g(\gamma)|_{g} \sim 1$ as $t \rightarrow 0$. Hence, 
$$(\langle \dot g(\gamma), \nabla_{N}\kappa \rangle_{g}  -  
\langle \nabla_{N}\dot g(\gamma), \kappa \rangle)_{g}dV_{S(t)} = 
\langle \nabla_{N}\kappa , \dot \gamma \rangle_{\gamma}dV_{S(t)} 
+ O(t)$$
$$= \langle \nabla_{N}\widetilde \kappa - n\widetilde \kappa, \dot \gamma 
\rangle_{\gamma} dV_{\gamma} + O(t).$$
Thus, 
\begin{equation}\label{e5.34}
\int_{B(t)}\langle (D_{1}\Phi)_{g}(\dot g(\gamma), \kappa \rangle dV_{g} = 
{\tfrac{1}{2}}\int_{S(t)}\langle \nabla_{N}\widetilde \kappa - n\widetilde \kappa, 
\dot \gamma \rangle_{\gamma}dV_{\gamma} + O(t).
\end{equation}
Now suppose, (contrary to (5.31)), 
\begin{equation}\label{e5.35}
\nabla_{N}\widetilde \kappa - n\widetilde \kappa = O(t),
\end{equation}
as forms on $(S(t), \bar g)$; note however that $\nabla$ is taken with respect to $g$ 
in (5.35). Since $\widetilde \kappa(N) = O(t)$ by (5.25), a simple computation shows that 
(5.35) implies that $\frac{1}{2}N(|\widetilde \kappa|^{2}) - n|\widetilde \kappa|^{2} = 
O(t)$. Integrating with respect to the induced metric on $(S(t), \bar g)$ and using the 
fact that $\frac{d}{dt}dV_{S(t)} = O(t)$, it follows that 
\begin{equation}\label{e5.36}
{\tfrac{1}{2}}N\int_{S(t)}|\widetilde \kappa|^{2}dV_{\gamma} - 
n\int_{S(t)}|\widetilde \kappa|^{2}dV_{\gamma} = O(t),
\end{equation}
as $t \rightarrow 0$. An elementary integration of (5.36) in $t$ then implies
$$\int_{S(t)}|\widetilde \kappa|^{2}dV_{\gamma} = O(t),$$
and hence, using (5.25) again, 
\begin{equation}\label{e5.37}
\kappa = o(t^{n}).
\end{equation}
(More precisely (5.37) holds in $L^{2}$, but pointwise decay then follows from elliptic 
regularity in weighted H\"older spaces, as in [26] or [28] for instance). Now $\kappa$ 
is an infinitesimal Einstein deformation, divergence-free by (5.26). Then Corollary 4.4, 
(cf.~also [5, Thm.~3.1]), and (5.37), together with the assumption in Theorem 1.3 that 
$\pi_{1}(M, \partial M) = 0$ imply that 
$$\kappa = 0 \ \ {\rm on} \ \ M,$$
giving a contradiction. This proves the relation (5.31). 

  The proof above shows that the form
\begin{equation}\label{e5.38}
\dot g(\gamma) = \lim_{t\rightarrow 0}\widetilde \kappa|_{S(t)},
\end{equation}
satisfies (5.31). The limit here exists by the regularity results of [16]. 
Thus, the space
\begin{equation}\label{e5.39}
\widetilde K = \{\widetilde \kappa = \lim_{t\rightarrow 0}t^{-n}\kappa|_{S(t)}: 
\kappa \in K\},
\end{equation}
is naturally identified with the cokernel of $D\Pi_{g}$ in $T_{\gamma}Met(\partial M)$. 
Note that $dim \widetilde K = dim K$ and also that the estimates (5.25) show that 
$\widetilde \kappa = \widetilde \kappa^{T}$ on $\partial M$. 

  This means that infinitesimal deformations of the boundary metric $\gamma$ in the 
direction $\widetilde \kappa$, $\widetilde \kappa \in \widetilde K$, are 
not realized as $\frac{d}{ds}\Pi(g_{s})|_{s=0}$, where $g_{s}$ is a curve in $E_{AH}$ 
through $g$, i.e.~a curve of {\it global} Einstein metrics on $M$. However, it is 
easy to see that $\widetilde \kappa$ is realized as the boundary variation of 
{\it locally defined} Einstein metrics. More precisely, choose any sequence 
$t_{i} \rightarrow 0$ and consider the metrics
\begin{equation}\label{e5.40}
g_{s,i} = g + st_{i}^{-n}\kappa + O(s^{2}),
\end{equation}
in the region $A_{i} = A(\frac{t_{i}}{2}, 2t_{i}) = \{x \in M: \frac{t_{i}}{2} \leq t(x) 
\leq 2t_{i}\}$, with $\kappa \in K$. For each $i$, this is a curve of metrics on $A_{i}$, 
Einstein to $1^{\rm st}$ order in $s$ at $s = 0$. The induced variation of the boundary 
metric on $S(t_{i})$ is, by construction, $(\widetilde \kappa)^{T}|_{S(t_{i})} \sim 
\widetilde \kappa|_{S(t_{i})}$. Now note that the linearized divergence constraint 
(5.6) or (5.8) only involves the behavior at $\partial M$, or equivalently, the limiting 
behavior on $(S(t_{i}), g_{t_{i}})$, $g_{t_{i}} = g|_{S(t_{i})}$, as $t_{i} \rightarrow 0$. 
This shows that the constraint (5.8) may be solved on the cokernel $\widetilde K$, and 
hence is solvable on all of $T_{\gamma}Met(\partial M)$. This proves the following:

\begin{corollary}\label{c5.5}
Let $g$ be a conformally compact Einstein metric on a compact manifold $M$ with 
$C^{\infty}$ boundary metric $\gamma$. Then the linearized divergence constraint 
equation {\rm (5.8)} is always solvable on $(\partial M, \gamma)$, i.e.~the map $\pi$ 
in {\rm (5.5)} is locally surjective at $(\gamma, \tau_{(n)})$. 
\end{corollary}

{\endproof}

\noindent
(A more detailed and computational proof of Corollary 5.5 will also be given 
below). 

  Combining Proposition 5.4 and Corollary 5.5, it follows that
\begin{equation}\label{e5.41}
{\mathcal L}_{X}\tau_{(n)} = 0, \ {\rm and \ hence} \ {\mathcal L}_{X}g_{(n)} = 0
\end{equation}
on $(\partial M, \gamma)$. (The second statement follows since the term $r_{(n)}$ is 
intrinsic to the boundary metric $\gamma$, so that ${\mathcal L}_{X}r_{(n)} = 0$). 
Proposition 5.1 and the assumption $\pi_{1}(M, \partial M) = 0$ then implies that 
$X$ extends to a Killing field $Y$ on $M$. This completes the first proof of 
Theorem 1.3. 

   Regarding the claim (5.29), since $\delta^{*}Y = 0$ and $Y$ is asymptotic to 
$X$ on $\partial M$, one has $\kappa = \delta^{*}(X - Y)$ on $M$, which implies 
that $\kappa = o(t^{n})$. By Corollary 4.4, this implies that $\kappa = 0$, 
as claimed. 

\medskip

  It is useful and of interest to give another, direct computational proof of 
Theorem 1.3, without using the identification (5.39) as the cokernel of $D\Pi$. 
The basic idea is to compute as in Proposition 5.4 on $(S(t), g_{t})$, with 
$A - Hg_{t}$ in place of $\tau_{(n)}$, and then pass to the limit on $\partial M$. 
Throughout the proof, we assume (5.28) holds. 

   Before starting the proof per se, we note that the estimates (5.25) and (5.28) 
imply that $X$ is tangential, i.e.~tangential to $(S(t), g)$, to high order, in 
that
\begin{equation}\label{e5.42}
\langle X, N \rangle = O(t^{n+1+\mu}).
\end{equation}
To see this, one has $(\delta^{*}X)(N, N) = \langle \nabla_{N}X, N \rangle = 
N \langle X, N \rangle$. Thus (5.42) follows from (5.25) and the claim that 
$\langle X, N \rangle = 0$ on $\partial M$. To prove the latter, consider the 
compactified metric $\bar g = t^{2}g$. One has ${\mathcal L}_{X}\bar g = 
{\mathcal L}_{X}(t^{2}g) = 2\frac{X(t)}{t}\bar g + O(t^{n})$. Thus for the 
induced metric $\gamma$ on $\partial M$, ${\mathcal L}_{X}\gamma = 
2\lambda \gamma$, where $\lambda = \lim_{t\rightarrow 0}\frac{X(t)}{t}$. 
Since $X$ is a Killing field on $(\partial M, \gamma)$, this gives $\lambda 
= 0$, which is equivalent to the statement that $\lim_{t\rightarrow 0}
\langle X, N \rangle_{g} = 0$. Note also that since $X$ is smooth up to 
$\partial M$, $|X|_{g} = O(t^{-1})$. 

  We claim also that
\begin{equation}\label{e5.43}
[X, N] = O(t^{n+1}),
\end{equation}
in norm. First, $\langle [X, N], N\rangle = \langle \nabla_{X}N - \nabla_{N}X, N \rangle 
= - (\delta^{*}X)(N,N) = O(t^{n+1+\mu})$. On the other hand, on tangential $g$-unit vectors 
$Y$, $\langle [X, N], Y\rangle = \langle \nabla_{X}N - \nabla_{N}X, Y \rangle 
\sim \langle \nabla_{X}N, Y \rangle - 2(\delta^{*}X)(N,Y) + \langle \nabla_{Y}X, N \rangle 
\sim - 2(\delta^{*}X)(N,Y) =  O(t^{n+1})$, as claimed. Here $\sim$ denotes equality modulo 
terms of order $o(t^{n})$. We have also used the fact that $\langle \nabla_{X}N, Y \rangle 
+ \langle \nabla_{Y}X, N \rangle \sim X\langle N, Y \rangle = 0$. 

   Now, to begin the proof itself, (assuming (5.28)), as above write
$$g_{s} = g + s\kappa + O(s^{2}) = g + s\delta^{*}X + O(s^{2}).$$
If $t_{s}$ is the geodesic defining function for $g_{s}$, (with boundary metric 
$\gamma$), then the Fefferman-Graham expansion gives $\bar g_{s} = dt_{s}^{2} + 
(\gamma + t_{s}^{2}g_{(2),s} + \dots + t_{s}^{n}g_{(n),s}) + O(t^{n+1})$. The 
estimate (5.42) implies that $t_{s} = t + sO(t^{n+2+\alpha}) + O(s^{2})$, 
so that modulo lower order terms, we may view $t_{s} \sim t$. Taking the derivative 
of the FG expansion with respect to $s$ at $s = 0$, and using the fact that 
$X$ is Killing on $(\partial M, \gamma)$, together with the fact that the lower 
order terms $g_{(k)}$, $k < n$, are determined by $\gamma$, it follows that, for 
$\widetilde k$ as in (5.39), 
\begin{equation}\label{e5.44}
\widetilde \kappa = {\tfrac{1}{2}}{\mathcal L}_{X}g_{(n)},
\end{equation}
at $\partial M$. Here both $\widetilde \kappa$ and ${\mathcal L}_{X}g_{(n)}$ are 
viewed as forms on $(\partial M, \gamma)$. 

   Next, we claim that on $(S(t), g_{t})$, 
\begin{equation}\label{e5.45}
{\mathcal L}_{X}A = -{\tfrac{n-2}{2}}t^{n-2}{\mathcal L}_{X}g_{(n)} + O(t^{n-1}),
\end{equation}
To see this, one has $A = \frac{1}{2}{\mathcal L}_{N}g = 
-\frac{1}{2}{\mathcal L}_{t\partial t}g = -\frac{1}{2}{\mathcal L}_{t\partial t}(t^{-2}g_{t})$. 
But ${\mathcal L}_{t\partial t}(t^{-2}g_{t}) = \sum {\mathcal L}_{t\partial t}(t^{-2+k}g_{(k)}) = 
\sum (k-2)t^{k-2}g_{(k)}$. The same reasoning as before then gives (5.45). 

  Given these results, we now compute
$$\int_{S(t)}\langle {\mathcal L}_{X}(A - Hg_{t}), \widetilde \kappa \rangle_{g_{t}} dV_{S(t)};$$
compare with the left side of (5.9). First, by (5.45),
$$\int_{S(t)}\langle {\mathcal L}_{X}A , \widetilde \kappa \rangle_{g_{t}} dV_{S(t)} = 
-{\tfrac{n-2}{2}}\int_{S(t)}\langle {\mathcal L}_{X}g_{(n)}, \widetilde \kappa \rangle_{\gamma} 
dV_{\gamma} + O(t).$$
Next, one has ${\mathcal L}_{X}(Hg_{t}) = X(H)g_{t} + H{\mathcal L}_{X}g_{t}$. For the 
first term, $X(H) = tr {\mathcal L}_{X}A + O(t^{n}) = -\frac{n-2}{2}t^{n-2}tr {\mathcal L}_{X}g_{(n)} 
+ O(t^{n})$. Since $tr g_{(n)}$ is intrinsic to $\gamma$ and $X$ is Killing on 
$(\partial M, \gamma)$, it follows that $X(H) = O(t^{n-1})$. Also, $\langle g_{t}, 
\widetilde \kappa \rangle = tr^{T}\widetilde \kappa$, where $tr^{T}$ is the tangential 
trace. By (5.29) and the fact that $\kappa$ is trace-free, $\langle g_{t}, \widetilde \kappa 
\rangle = O(t^{1+\alpha})$. Hence $X(H)\langle g_{t}, \widetilde \kappa \rangle dV_{S(t)} = 
O(t^{\alpha})$. Similarly, from (5.43) one computes ${\mathcal L}_{X}g_{t} = {\mathcal L}_{X}g 
+ O(t^{n+1}) = 2t^{n}\widetilde \kappa + O(t^{n+1})$. Since $H \sim n$, using (5.44) this gives
$$-\int_{S(t)}\langle {\mathcal L}_{X}(Hg_{t}), \widetilde \kappa \rangle dV_{S(t)} = 
-n\int_{S(t)}\langle {\mathcal L}_{X}g_{(n)}, \widetilde \kappa \rangle_{\gamma} 
dV_{\gamma} + O(t^{\alpha}).$$
Combining these computations then gives
\begin{equation}\label{e5.46}
\int_{S(t)}\langle {\mathcal L}_{X}(A - Hg_{t}), \widetilde \kappa \rangle_{g_{t}} dV_{S(t)} = 
-({\tfrac{n-2}{2}}+n)\int_{\partial M}\langle {\mathcal L}_{X}g_{(n)}, 
\widetilde \kappa \rangle_{\gamma} dV_{\gamma} + o(1).
\end{equation}

  On the other hand, one may use (5.12), with $h = \widetilde \kappa$, to compute the left 
side of (5.46). Since $\kappa$ is transverse-traceless, the estimate (5.42) implies 
that $tr^{T}\widetilde \kappa = O(t^{n+1})$ and $\delta^{T} X = O(t^{n+1})$. Hence the 
last two terms in (5.12) are lower order, (i.e.~$O(t)$). Also, as in (5.40), the metrics 
$g_{s}$ are Einstein to $1^{\rm st}$ order in $s$ in a neighborhood of $S(t)$ and hence
\begin{equation}\label{e5.47}
\int_{S(t)}\langle {\mathcal L}_{X}(A - Hg_{t}), \widetilde \kappa \rangle_{g_{t}} dV_{S(t)} = 
-2\int_{S(t)}\langle \delta'(A - Hg_{t}), X \rangle dV_{g_{t}} = 
\end{equation}
$$2\int_{S(t)}\langle (A - Hg_{t})', (\delta^{*})^{T}X \rangle dV_{g_{t}} + O(t).$$
Now by (5.42), $(\delta^{*})^{T}X = \delta^{*}X + O(t^{n+1})$ while $A' = 
\frac{d}{ds}(A_{g+s\widetilde \kappa}) = \frac{1}{2}({\mathcal L}_{N}\widetilde \kappa + 
{\mathcal L}_{N'}g) = \frac{1}{2}\nabla_{N}\widetilde \kappa + O(t)$. It follows that
\begin{equation}\label{e5.48}
\int_{S(t)}\langle (A - Hg_{t})', (\delta^{*})^{T}X \rangle dV_{g_{t}} = 
{\tfrac{1}{2}}\int_{S(t)}\langle \nabla_{N}\widetilde \kappa, \widetilde \kappa 
\rangle_{\bar g_{t}} dV_{\gamma} + O(t) = O(t),
\end{equation}
where the last equality follows by the same reasoning as following (5.36), (since 
$N = -t\partial_{t}$). 

  The equations (5.46)-(5.48) imply that
$$\int_{S(t)}\langle {\mathcal L}_{X}g_{(n)}, \widetilde \kappa \rangle_{\gamma} dV_{\gamma} = 
O(t).$$
Via (5.44) again, this shows that 
$$\widetilde \kappa = 0$$
on $\partial M$, and hence via Corollary 4.4 again, $\kappa = 0$ on $M$. This completes 
the second proof of Theorem 1.3. 
{\endproof}

 {\bf Proof of Corollary 1.4.}

   Suppose $(M, g)$ is a conformally compact Einstein metric with boundary metric 
given by the round metric $S^{n}(1)$ on $S^{n}$. Theorem 1.3 implies that the 
isometry group of $(M, g)$ contains the isometry group of $S^{n}$. This reduces 
the Einstein equations to a simple system of ODE's, and it is easily seen that 
the only solution is given by the Poincar\'e metric on the ball $B^{n+1}$. 
{\endproof}

\begin{remark}\label{r5.6}
{\rm By means of Obata's theorem [30], Theorem 1.3 remains true for continuous 
groups of conformal isometries at conformal infinity. Thus, the class of the 
round metric on $S^{n}$ is the only conformal class which supports a non-essential 
conformal Killing field, i.e.~a field which is not Killing with respect to some 
conformally related metric. Corollary 1.4 shows that any $g \in E_{AH}$ with 
boundary metric $S^{n}(1)$ is necessarily the hyperbolic metric $g_{-1}$ on 
the ball. For $g_{-1}$, it is well-known that essential conformal Killing fields 
on $S^{n}$ extend to Killing fields on $({\mathbb H}^{n+1}, g_{-1})$. 

  We expect that a modification of the proof of Theorem 1.3 would give this result 
directly, without the use of Obata's theorem. In fact, such would probably give 
(yet) another proof of Obata's result. }
\end{remark}

  Corollary 5.5 shows, in the global situation, that the projection $\pi$ of 
the constraint manifold ${\mathcal T}$ to $Met(\partial M)$ is always locally surjective. 
Hence there exists a formal solution, and an exact solution in the analytic case, for 
any nearby boundary metric, which is defined in a neighborhood of the boundary. However, 
the full boundary map $\Pi$ in (5.13) or (5.14) on global metrics is not locally surjective 
in general; nor is it always globally surjective. 

  The simplest example of this behavior is provided by the family of AdS Schwarzschild 
metrics. These are metrics on ${\mathbb R}^{2}\times S^{n-1}$ of the form
$$g_{m} = V^{-1}dr^{2} + Vd\theta^{2} + r^{2}g_{S^{n-1}(1)},$$
where $V = V(r) = 1 + r^{2} - \frac{2m}{r^{n-2}}$. Here $m > 0$ and $r \in [r_{+}, \infty]$, 
where $r_{+}$ is the largest root of the equation $V(r_{+}) = 0$. The locus $\{r_{+} = 0\}$ 
is a totally geodesic round $S^{n-1}$ of radius $r_{+}$. Smoothness of the metric at 
$\{r_{+} = 0\}$ requires that the circular parameter $\theta$ runs over the interval 
$[0,\beta]$, where 
$$\beta = \frac{4\pi r_{+}}{nr_{+}^{2}+(n-2)}.$$
The metrics $g_{m}$ are isometrically distinct for distinct values of $m$, and form 
a curve in $E_{AH}$ with conformal infinity given by the conformal class of the 
product metric on $S^{1}(\beta)\times S^{n-1}(1)$. As $m$ ranges over the interval 
$(0, \infty)$, $\beta$ has a maximum value of 
$$\beta \leq \beta_{\max} = 2\pi \sqrt{(n-2)/n}.$$
As $m \rightarrow 0$ or $m \rightarrow \infty$, $\beta \rightarrow 0$. 

   Hence, the metrics $S^{1}(L)\times S^{n-1}(1)$ are not in $\Pi(g_{m})$ for any 
$L > \beta_{\max}$. In fact these boundary metrics are not in $Im(\Pi)$ generally, 
for any manifold $M^{n+1}$. For Theorem 1.3 implies that any conformally compact 
Einstein metric with boundary metric $S^{1}(L)\times S^{n-1}(1)$ has an isometry 
group containing the isometry group of $S^{1}(L)\times S^{n-1}(1)$. This again 
reduces the Einstein equations to a system of ODE's and it is easy to see, (although 
we do not give the calculations here), that any such metric is an AdS Schwarzschild metric.

\begin{remark}\label{r5.7}
{\rm In the context of Propositions 5.1 and 5.4, it is natural to consider the 
issue of whether local Killing fields of $\partial M$, (i.e.~Killing fields defined 
on the universal cover), extend to local Killing fields of any global conformally 
compact Einstein metric. Note that Proposition 5.1 and Proposition 5.4 are both 
local results, the latter by using variations $h_{(0)}$ which are of compact support. 
However, the linearized constraint condition (5.8) is not invariant under covering 
spaces; even the splitting (5.7) is not invariant under coverings, since a Killing 
field on a covering space need not descend to the base space. 

  We claim that local Killing fields do not extend even locally into the interior in 
general. As a specific example, let $N^{n+1}$ be any complete, geometrically finite 
hyperbolic manifold, with conformal infinity $(\partial N, \gamma)$, and which has 
at least one parabolic end, i.e.~a finite volume cusp end, with cross sections given 
by flat tori $T^{n}$. There exist many such manifolds. The metric at conformal infinity 
is conformally flat, so there are many local Killing fields on $\partial N$. For example, 
in many examples $N$ itself is a compact hyperbolic manifold. Of course the local 
(conformal) isometries of $\partial N$ extend here to local isometries of $N$. 

  However, as shown in [15], the cusp end may be capped off by Dehn filling with a 
solid torus, to give infinitely many distinct conformally compact Einstein metrics 
with the same boundary metric $(\partial N, \gamma)$. These Dehn-filled Einstein metrics 
cannot inherit all the local conformal symmetries of the boundary. }
\end{remark}

\begin{remark} \label{r5.8} 
{\rm We point out that Theorem 1.3 fails for complete Ricci-flat metrics which are ALE 
(asymptotically locally Euclidean). The simplest counterexamples are the family of 
Eguchi-Hanson metrics, which have boundary metric at infinity given by the round metric 
on $S^{3}/{\mathbb Z}_{2}$. The symmetry group of these metrics is strictly smaller than 
the isometry group $Isom (S^{3}/{\mathbb Z}_{2})$ of the boundary. Similarly, the 
Gibbons-Hawking family of metrics with boundary metric the round metric on 
$S^{3}/{\mathbb Z}_{k}$ have only an $S^{1}$ isometry group, much smaller than 
the group $Isom (S^{3}/{\mathbb Z}_{k})$. 

  This indicates that, despite a number of proposals, some important features of holographic 
renormalization in the AdS context cannot carry over to the asymptotically flat case. }
\end{remark}

\bibliographystyle{plain}

\begin{center}
March, 2006/May, 2007
\end{center}

\noindent
{\address Department of Mathematics\\
S.U.N.Y. at Stony Brook\\
Stony Brook, NY 11794-3651\\
E-mail: anderson@math.sunysb.edu}

\end{document}